# Multicomputation with Numbers: The Case of Simple Multiway Systems

Stephen Wolfram*

*Integer iteration rules such as $n \longmapsto \{a\,n + b, c\,n + d\}$ are studied as minimal examples of the general process of multicomputation. Despite the simplicity of such rules, their multiway graphs can be complex, exhibiting, for example, emergent geometry and difficult questions of confluence. Generalizations to rules involving non-integers and other functions are also considered. Connections with physics and with various number-theoretic and other questions are made.*

## A Minimal Example of Multicomputation

Multicomputation is an important new paradigm, but one that can be quite difficult to understand. Here my goal is to discuss a minimal example: multiway systems based on numbers. Many general multicomputational phenomena will show up here in simple forms (though others will not). And the involvement of numbers will often allow us to make immediate use of traditional mathematical methods.

A multiway system can be described as taking each of its states and repeatedly replacing it according to some rule or rules with a collection of states, merging any states produced that are identical. In our Physics Project, the states are combinations of relations between elements, represented by hypergraphs. We've also often considered string substitution systems, in which the states are strings of characters. But here I'll consider the case in which the states are numbers, and for now just single integers.

And in this case multiway systems can be represented in a particularly simple way, with each state $s$ just being repeatedly replaced according to:

$$s \rightarrow \{f_1(s), f_2(s), \ldots, f_p(s)\}$$

For a "binary branching" case the update rule is

$$s \longmapsto \{f_1[s],\ f_2[s]\}$$

and one can represent the evolution of the system by the multiway graph which begins:

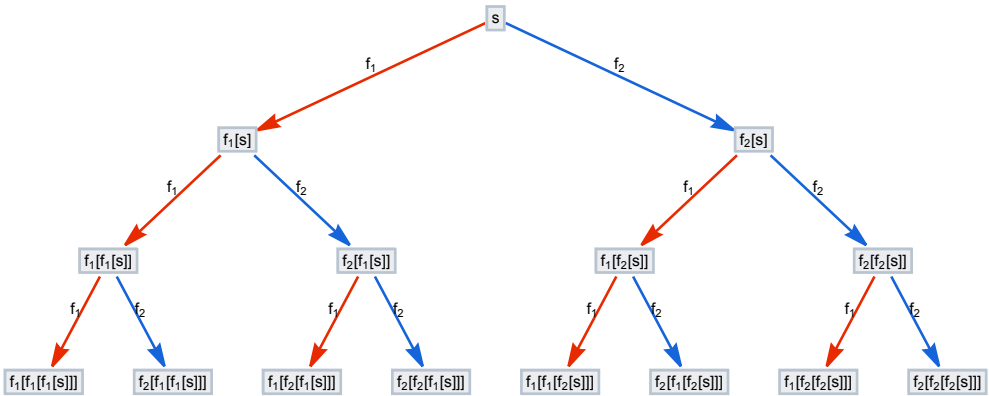

and continues (indicating $f_1$ by red and $f_2$ by blue):

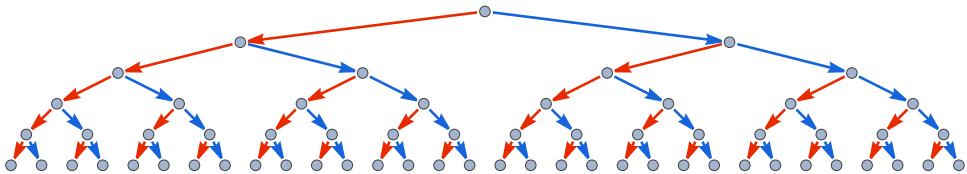

With arbitrary "symbolic" $f_i$ this ("free multiway system") tree is the only structure one can get. But things can get much less trivial when there are forms for $f_1$, $f_2$ that "evaluate" in some way, because then there can be identities that make branches merge. And indeed most of what we'll be discussing here is associated with this phenomenon and with the "entanglements" between states to which it leads.

It's worth noting that the specific setup we're using here avoids quite a lot of the structural complexity that can exist in multicomputational systems. In the general case, states can contain multiple "tokens", and updates can also "consume" multiple tokens. In our case here, each state just contains one token—which is a single number—and this is what is "consumed" at each step. (In our Physics Project, a state corresponds to a hyperedge which contains many hyperedge tokens, and the update rule typically consumes multiple hyperedges. In a string substitution system, a state is a character string which contains many character tokens, and the update typically consumes multiple—in this case, adjacent—character tokens.)

With the setup we're using here there's one input but multiple outputs (2 in the example above) each time the update rule is applied (with the inputs and outputs each being individual numbers). It's also perfectly possible to consider cases in which there are multiple inputs as well as multiple outputs. But here we'll restrict ourselves to the "one-to-many" ("traditional multiway") case. And it's notable that this case is exceptionally easy to describe in the Wolfram Language:

*In[ ]:=* NestGraph[s ⟼ {f₁[s], f₂[s]}, s, 3, VertexLabels → Automatic]

*Out[ ]=*

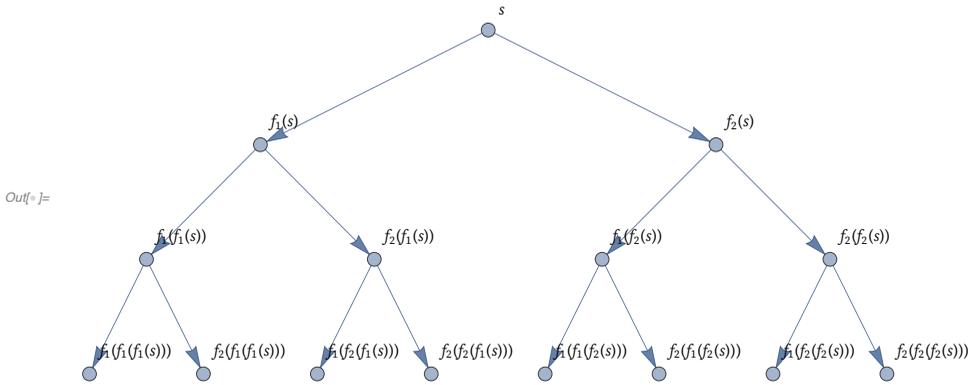

# Multiway Systems Based on Addition

As our first example, let's consider multiway systems whose rules just involve addition.

The trivial ("one-input, one-output") rule

$n \mapsto \{n + 1\}$

gives a multiway graph corresponding to a "one-way number line":

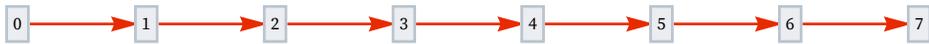

The rule

$n \mapsto \{n - 1, n + 1\}$

gives a "two-way number line":

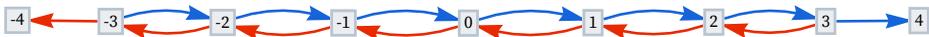

But even

$n \mapsto \{n + 1, n + 2\}$

gives a slightly more complicated multiway graph:

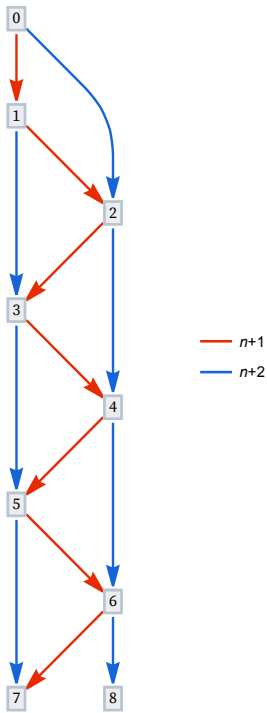

What's going on here? Basically each triangle represents an identity. For example, starting from 1, applying $n \mapsto n + 1$ twice gives 3, which is the same result as applying $n \mapsto n + 2$ once. Or, writing the rule in the form

$n \mapsto \{f_1[n], f_2[n]\}$

the triangles are all the result of the fact that in this case

$f_1[f_1[n]] = f_2[n]$

For the $n \mapsto \{n + 1\}$ "number line" rule, it's obvious that we'll eventually visit every integer—and the +1, +2 rule also visits every integer.

Consider now instead of +1 and +2 the case of +2 and +3:

$n \mapsto \{n + 2, n + 3\}$

After a few steps this gives:

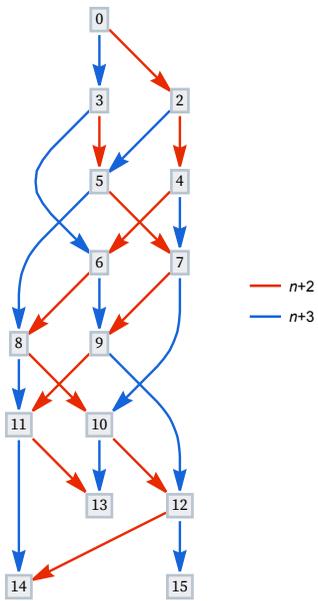

Continuing a little longer gives:

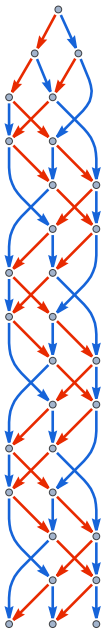

It's a little difficult to see what's going on here. It helps to show which edges correspond to +2 and +3:

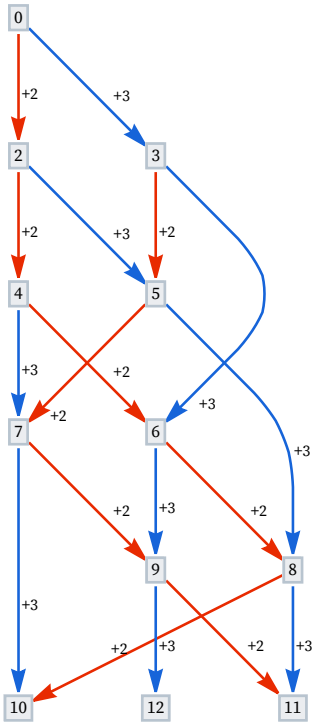

We'll return to this a little later, but once again we can see that there are cycles in this graph, corresponding to simple "commutativity identities", such as

$0 + 2 + 3 = 0 + 3 + 2$

and

$0 + 2 + 2 + 3 = 0 + 3 + 2 + 2$

as well as "LCM identities" such as

$0 + 2 + 2 + 2 = 0 + 3 + 3$

(Note that in this case, all integers above 1 are eventually generated.)

Let's look now at a case with slightly larger integers:

$n \longmapsto \{n + 4, \ n + 7\}$

After 6 steps one gets a simple grid

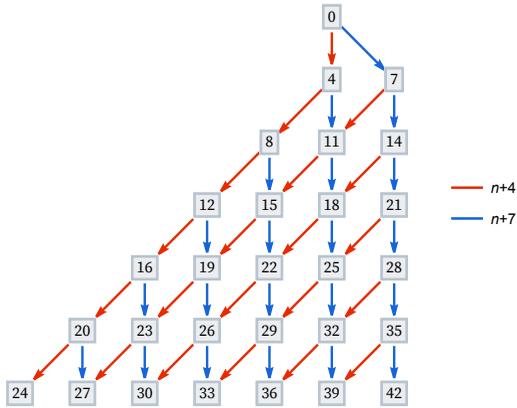

essentially made up of "commutativity identities". But continuing a little longer one sees that it begins to "wrap around"

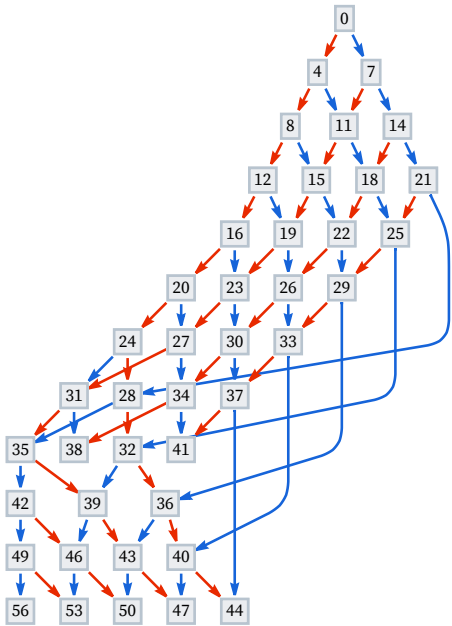

eventually forming a kind of "tube" with a spiral grid on the outside:

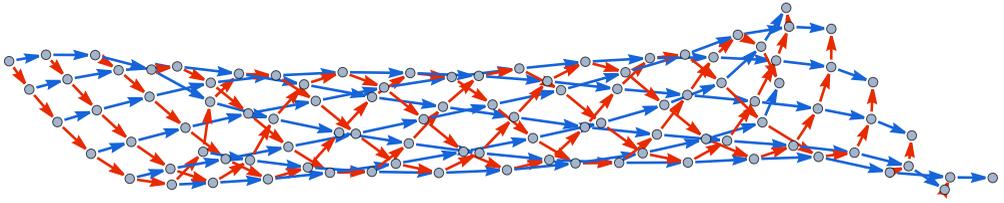

The "grid" is defined by "commutativity identities". But the reason it's a "closed tube" is that there are also "LCM identities". To understand this, unravel everything into a grid with +4 and +7 directions—then draw lines between the duplicated numbers:

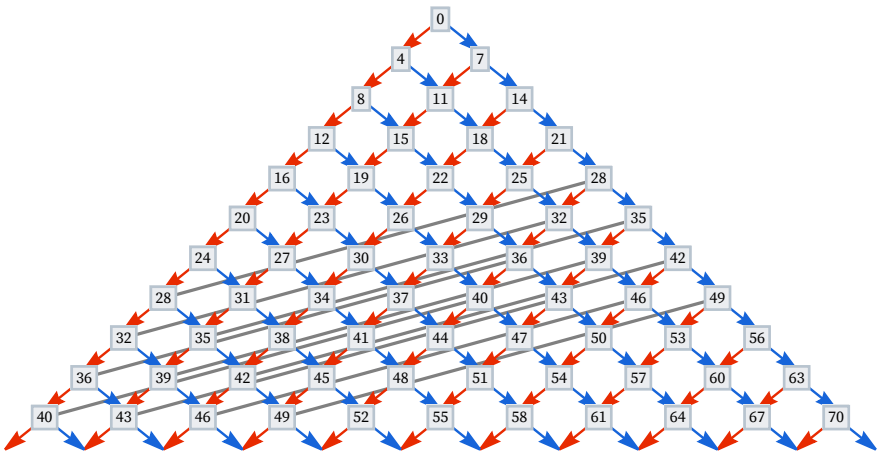

The "tube" is formed by rolling the grid up in such a way as to merge these numbers. But now if we assume that the multiway graph is laid out (in 3D) so that each graph edge has unit length, application of Pythagoras's theorem in the picture above shows that the effective circumference of the tube is $\sqrt{4^2 + 7^2} \approx 8.06$.

In another representation, we can unravel the tube by plotting numbers at {x, y} according to their decomposition in the form 7 x + y:

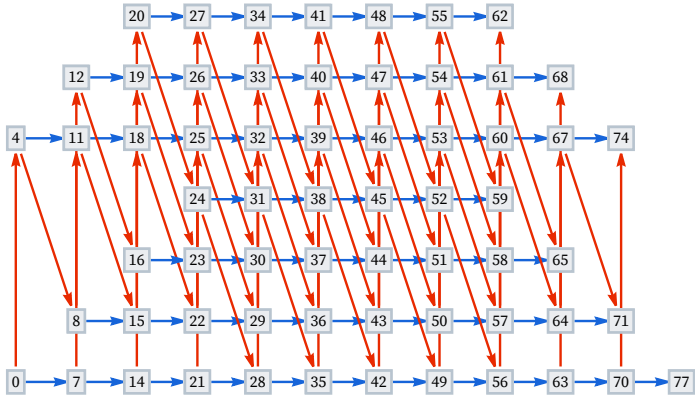

(From this representation we can see that every value of *n* can be reached so long as *n*>17.)

For the rule

$n \longmapsto \{n + 7, n + 11\}$

the multiway graph forms a tube of circumference $\sqrt{7^2+11^2} \approx 13$ which can be visualized in 3D as:

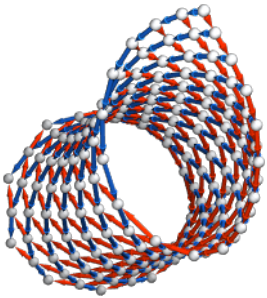

And what's notable here is that even though we're just following a simple discrete arithmetic process, we're somehow "inevitably getting geometry" out of it. It's a tiny, toy example of a much more general and powerful phenomenon that seems to be ubiquitous in multicomputational systems—and that in our models of physics is basically what leads to the emergence of things like the limiting continuum structure of space.

We've seen a few specific example of "multiway addition systems". What about the more general case?

For

$n \longmapsto \{n + a, n + b\}$

a "tube" is generated with circumference

$\sqrt{\bar{a}^2 + \bar{b}^2}$ where $\{\bar{a}, \bar{b}\} = \{a, b\} \, / \, \mathsf{GCD}[a, b]$

After enough steps, all integers of the form $k\,\mathsf{GCD}[a, b]$ will eventually be produced—which means that all integers are produced if $a$ and $b$ are relatively prime. There's always a threshold, however, given by $\mathsf{FrobeniusNumber}[\{a, b\}]$—which for $a$ and $b$ relatively prime is just $a\,b - a - b$.

By the way, a particular number $n$—if it's going to be generated at all—will first be generated at step

Min[Map[Total, FrobeniusSolve[{$a$, $b$}, $n$]]]

(Note that the fact that the multiway graph approximates a finite-radius tube is a consequence of the commensurability of any integers $a$ and $b$. If we had a rule like $n \longmapsto \{n + 1, n + \sqrt{2}\,\}$, we'd get an infinite 2D grid.)

For

$n \longmapsto \{n + a, \, n + b, \, n + c\}$

a tube is again formed, with a circumference effectively determined by the smaller pair (after $\mathsf{GCD}$ reduction) of $a$, $b$ and $c$. And if $\mathsf{GCD}[a, b, c] = 1$, all numbers above $\mathsf{FrobeniusNumber}[\{a, b, c\}]$ will eventually be generated.

## Pure Multiplication

One of the simplest cases of multiway systems are those based on pure multiplication. An example is (now starting from 1 rather than 0):

$n \longmapsto \{2\,n, \, 3\,n\}$

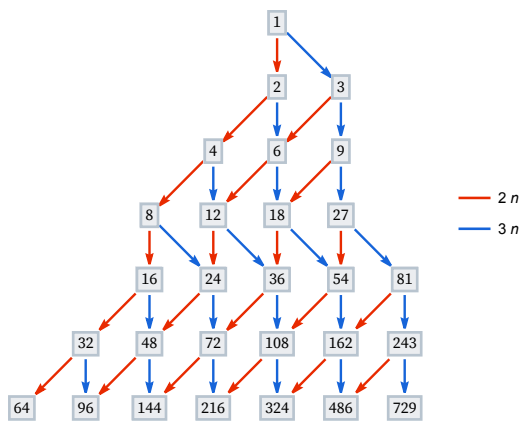

In general, for

$$n \longmapsto \{a\, n, b\, n\}$$

we'll get a simple 2D grid whenever $a$ and $b$ aren't both powers of the same number. With $d$ elements in the rule we'll get a $d$-dimensional grid. For example,

$$n \longmapsto \{2\, n, 3\, n, 5\, n\}$$

gives a 3D grid:

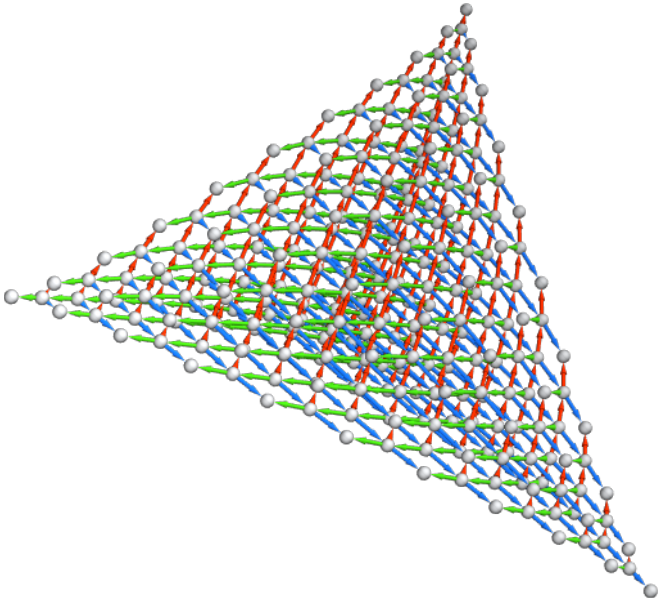

If the multipliers in the rule are all powers of the same number, the multiway graph degenerates to some kind of ladder. In the case

$$n \longmapsto \{n, 2\, n\}$$

this is just:

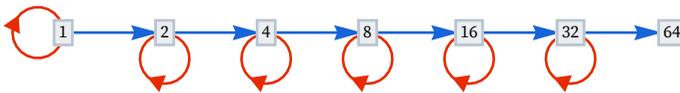

while for

$$n \longmapsto \{2\, n, 8\, n\}$$

it is

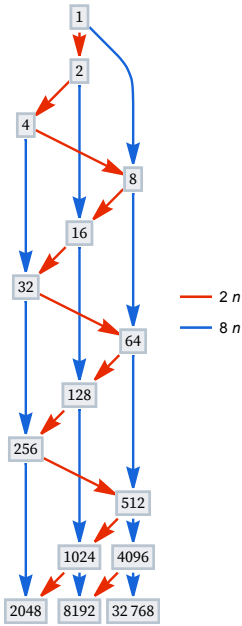

and in general for

$$n \longmapsto \{a\,n,\, a^m\,n\}$$

it is a "width-$m$" ladder graph.

## Multiplication and Addition: n ⟼ {a n, n + b}

Let's look now at combining multiplication and addition—to form what we might call affine multiway systems. As a first example, consider the case (which I actually already mentioned in *A New Kind of Science*):

$$n \longmapsto \{2\,n,\, n+1\}$$

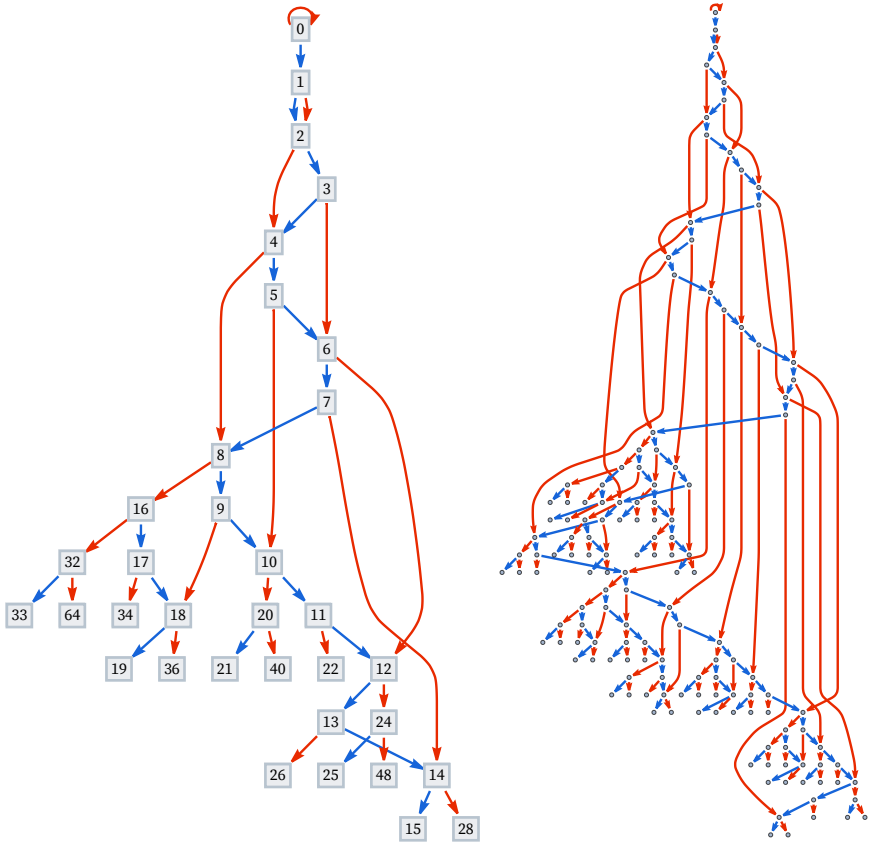

Considering the simplicity of the rule by which it was generated, this result looks surprisingly complex. One immediate result is that after $t$ steps, the total number of distinct numbers reached is **Fibonacci**[$t - 1$], which increases exponentially like $\phi^t$. Eventually the $n \mapsto n + 1$ ensures that every integer is generated. But the $n \mapsto 2\,n$ often "jumps ahead", and since the maximum number generated at step $t$ is $2^{t-2}$ the "average density" of numbers falls exponentially like $(\frac{\phi}{2})^t \approx 1.24^{-t}$.

Continuing the evolution further and using a different rendering we get the very "geometrical" (planar) structure

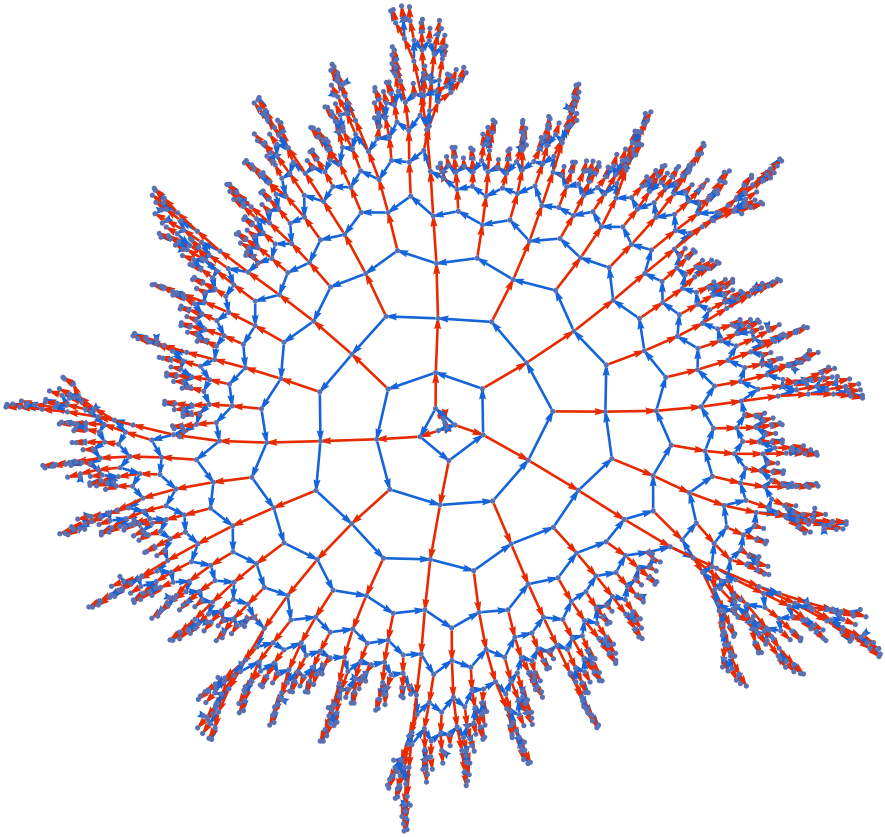

What can we say about this structure? Apart from the first few steps (rendered at the center), it consists of a spiral of pentagons. Each pentagon (except the one at the center) has the form

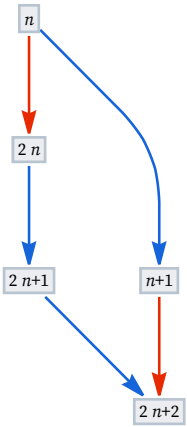

reflecting the relation

$$f_2[f_2[f_1[n]]] = f_1[f_2[n]]$$

Going out from the center, each successive layer in the spiral has twice the number of pentagons, with each pentagon at a given layer "spawning" two new pentagons at the next layer.

Removing "incomplete pentagons" this can be rendered as:

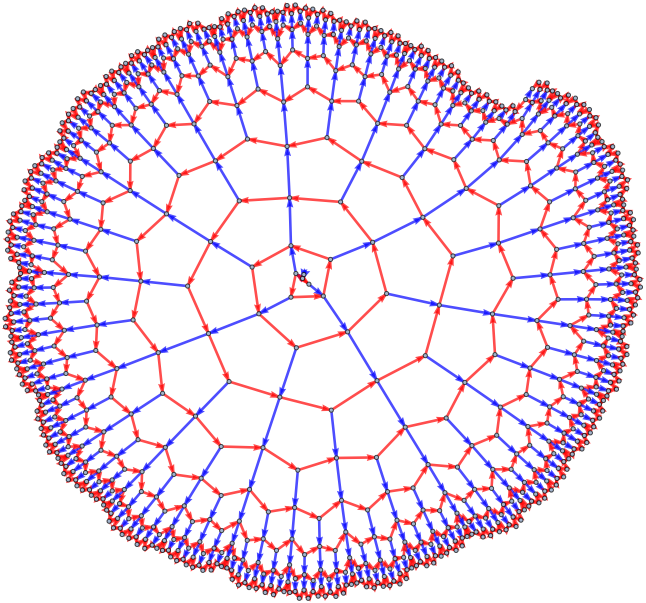

What about other rules of the general form:

$$n \mapsto \{a\,n,\ n+1\}$$

Here are the corresponding ("complete polygon") results for $a$=3 through 5:

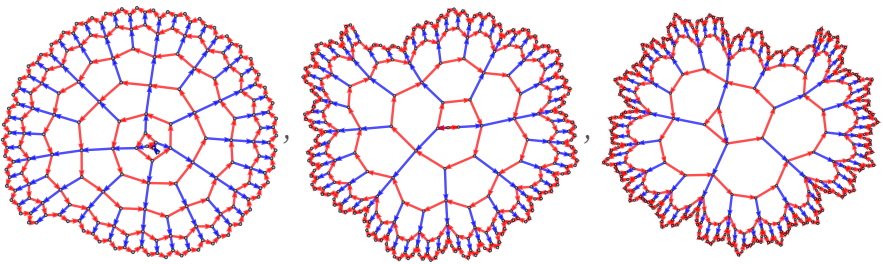

The multiway graphs in these cases correspond to spirals of $(a+3)$-gons defined by the identity

$$\mathsf{Nest}[f_2, f_1[n], a] = f_1[f_2[n]]$$

or equivalently

$$a\,n + a = (n+1)\,a$$

At successive layers in the spiral, the number of ($a$+3)-gons increases like $a^s$.

Eventually the evolution of the system generates all possible integers, but at step $t$ the number of distinct integers obtained so far is given by the generalized Fibonacci series obtained from

LinearRecurrence[Table[1, $a$], PadLeft[{1}, $a$], {$t$ + $a$}]

which for large $t$ is

$\phi_{a-1}{}^t$

where $\phi_k$ is the k-nacci generalized golden ratio, which approaches $2 - 2^{-k}$ for large $k$.

If we consider

$n \mapsto \{a\,n, \ n + b\}$

it turns out that one gets the same basic structure (with ($a$ + 3)-gons) for $b > 1$ as for $b = 1$. For example, with

$n \mapsto \{4\,n, \ n + 7\}$

one gets:

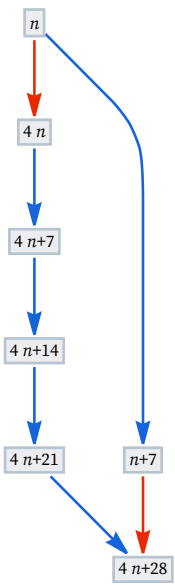

## The Rule n ⟼ {2n + 1, 3n + 1}

For the rule

$n \mapsto \{2\,n + 1, \ 3\,n + 1\}$

there are at first no equivalences that cause merging in the multiway graph:

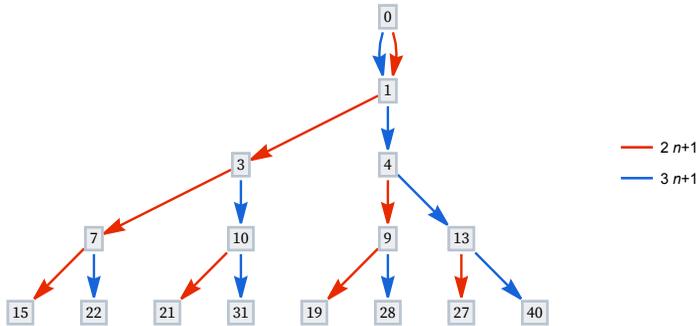

But after 5 steps we get

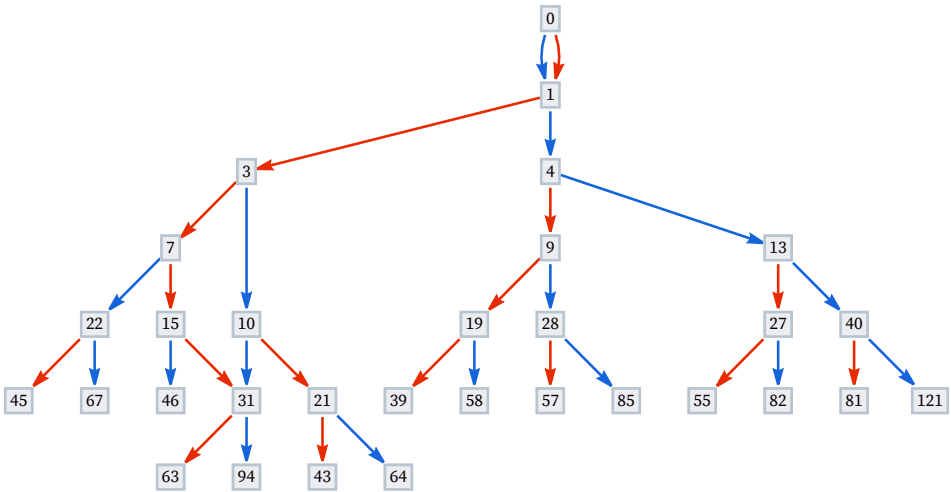

where now we see that 15 and 31 are connected "across branches".

After 10 steps this becomes:

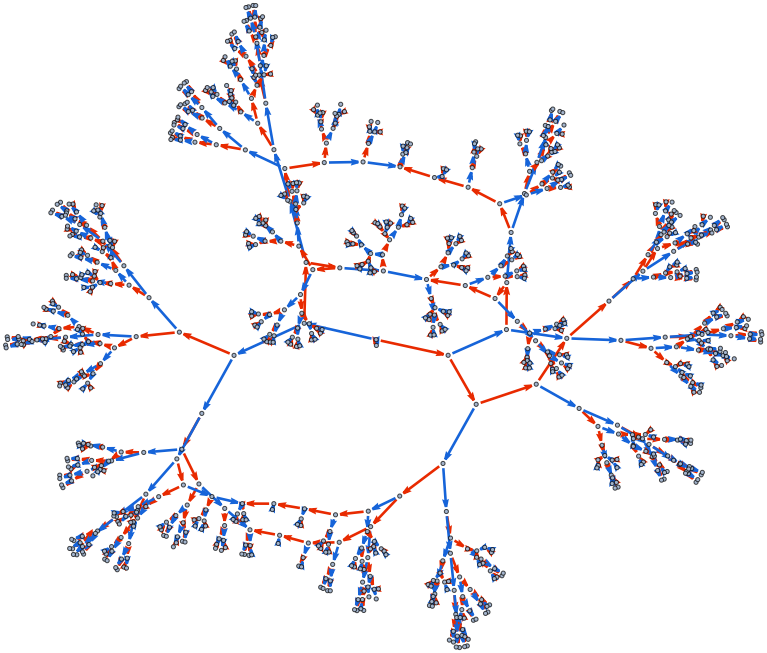

At a visual level this seems to consist of two basic components. First, a collection of loops, and second a collection of tree-like "loose ends". Keeping only complete loops and going a few more steps we get:

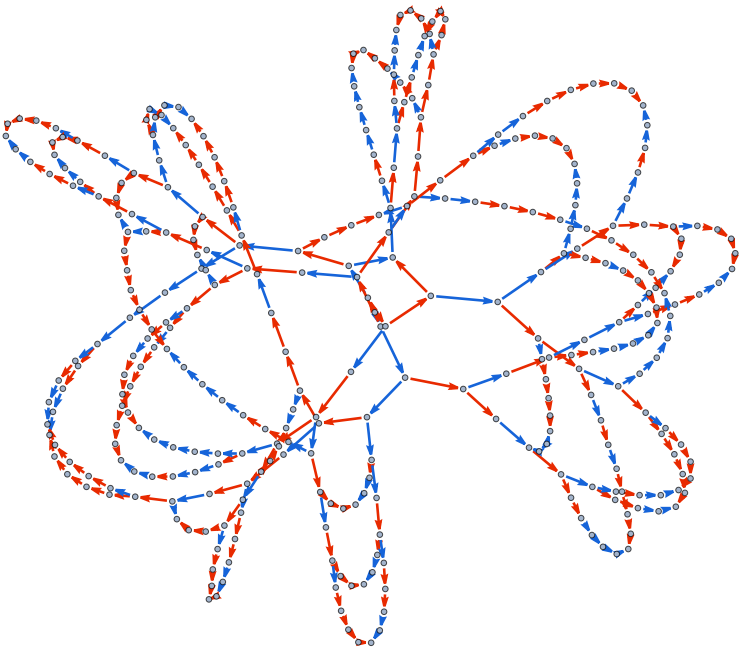

Unlike in previous cases, the "loops" (AKA "polygons") are not of constant size. Here are the first few that occur (note these loops "overlap" in the sense that several "start the same way"):

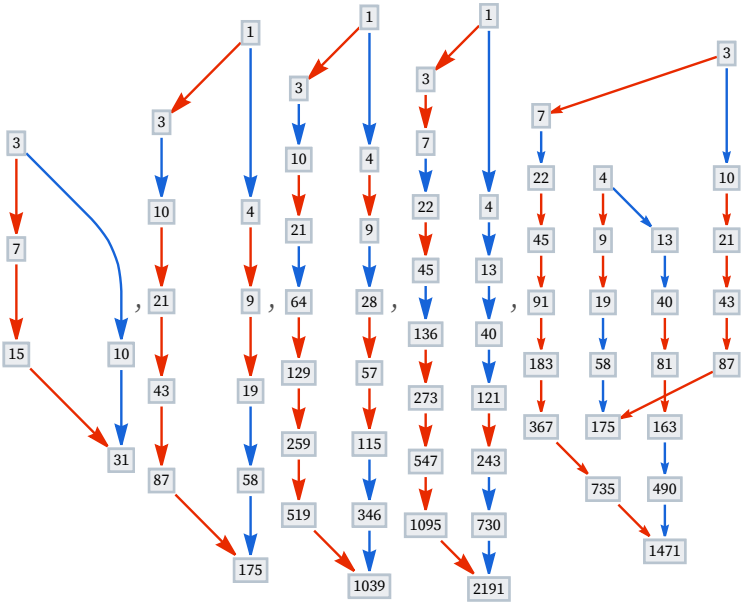

As before, each of these loops in effect corresponds to an identity about compositions of functions—though now it matters what these compositions are applied to. So, for example, the 4[th] loop above corresponds to (where $k$ stands for the function $n \mapsto k\,n + 1$):

223232222[1] == 3333233[1]

In explicit form this becomes:

2 (2 (2 (2 (2 (3 (2 (3 (2 (2 1+1)+1)+1)+1)+1)+1)+1)+1)+1)+1 = 3 (3 (2 (3 (3 (3 (3 1+1)+1)+1)+1)+1)+1)+1

where both sides evaluate to the same number, in this case 26815.

Much as in the Physics Project, we can think of each "loop" as beginning with the creation of a "branch pair", and ending with the merger of the different paths from each member of the pair. In a later section we'll discuss the question of whether every branch pair always in the end re-merges. But for now we can just enumerate mergers—and we find that the first few occur at:

{31, 175, 1039, 1471, 2191, 4495, 6223, 8815, 13 135, 20 479, 22 639, 26 815}

(Note that a merger can never involve more than two branches, since any given number has at most one "pre-image" under $n \mapsto 2\,n + 1$ and one under $n \mapsto 3\,n + 1$.)

Here is a plot of the positions of the $k^{\text{th}}$ mergers—together with a quadratic fit (indicated by the dotted line):

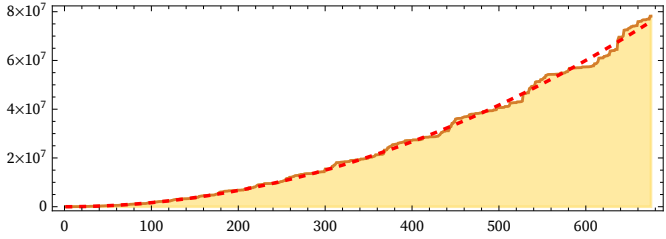

(As we'll discuss later, the numbers at which these mergers occur are for example always of the form $144\,u+31$.)

Taking second differences indicates a certain apparent randomness:

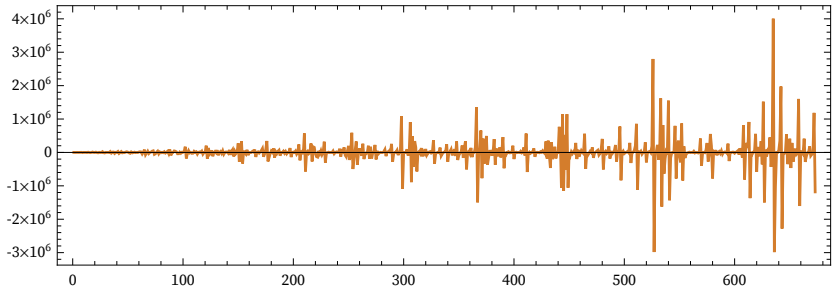

What can we say about the overall structure of the multiway graph? One basic question is what numbers ever even occur in the evolution of the system. Here are the first few, for evolution starting from 0:

{0, 1, 3, 4, 7, 9, 10, 13, 15, 19, 21, 22, 27, 28, 31, 39, 40, 43, 45, 46, 55, 57, 58, 63, 64}

And here are successive differences

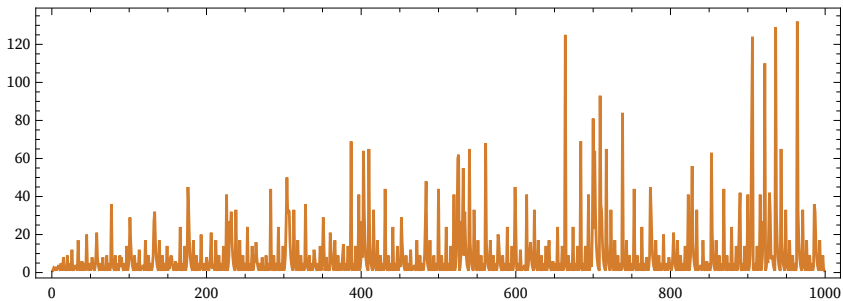

Dividing successive $m$ by the $m^{\text{th}}$ number gives a progressive estimate of the density of numbers:

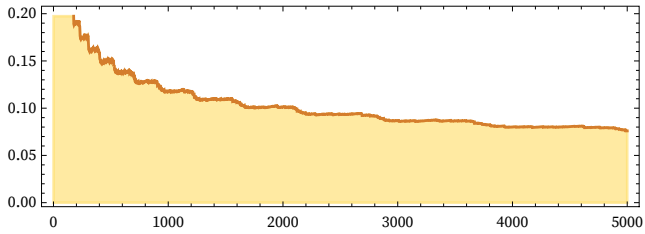

On a log-log scale this becomes

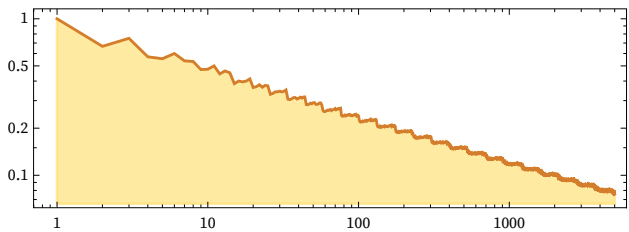

showing a rough fit to $m^{-0.3}$—and suggesting an asymptotic density of 0.

Note, by the way, that while the maximum gap grows on average linearly (roughly like $0.17\,m$)

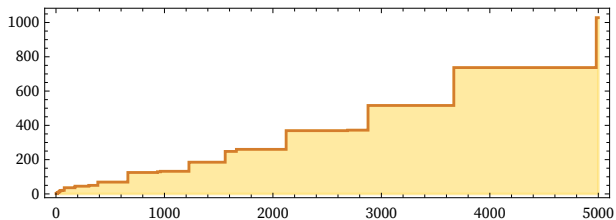

the distance between gaps of size 1 shows evidence of remaining bounded:

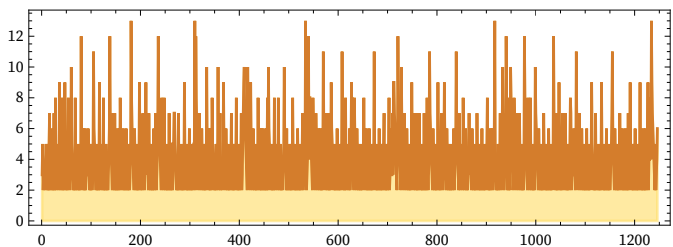

(A related result from the 1970s states that the original sequence contains infinite-length arithmetic progressions—implying the presence of infinite runs of numbers whose differences are constant.)

# The More General "Affine" Case: $n \longmapsto \{a\,n + b,\, c\,n + d\}$

Not every rule of the form

$n \longmapsto \{a\,n + b,\, c\,n + d\}$

leads to a complex multiway graph. For example

$n \longmapsto \{2\,n,\, 2\,n + 1\}$

just gives a pure binary tree since $2n$ just adds a 1 at the beginning of the binary digit sequence of $n$, while $2n{+}1$ adds one at the end:

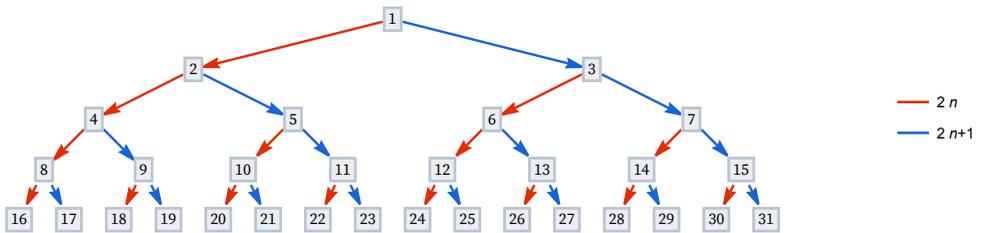

Meanwhile

$n \longmapsto \{2\,n + 1,\, 3\,n + 2\}$

gives a simple grid

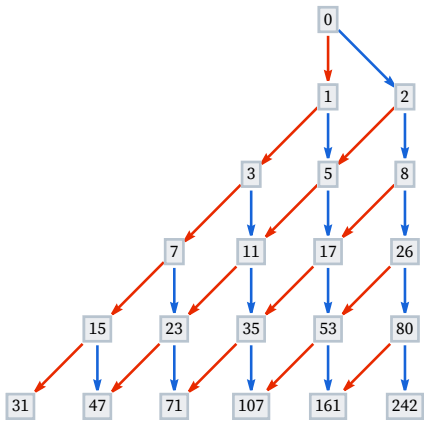

where at level $t$ the numbers that appear are simply

$-1 + 2^{t + 1 - i} \times 3^i$

and the pattern of use of the two cases in the rule makes it clear why the grid structure occurs.

Here are the behaviors of all inequivalent nontrivial rules of the form

$$n \mapsto \{a\,n + b,\ c\,n + d\}$$

with constants up to 3:

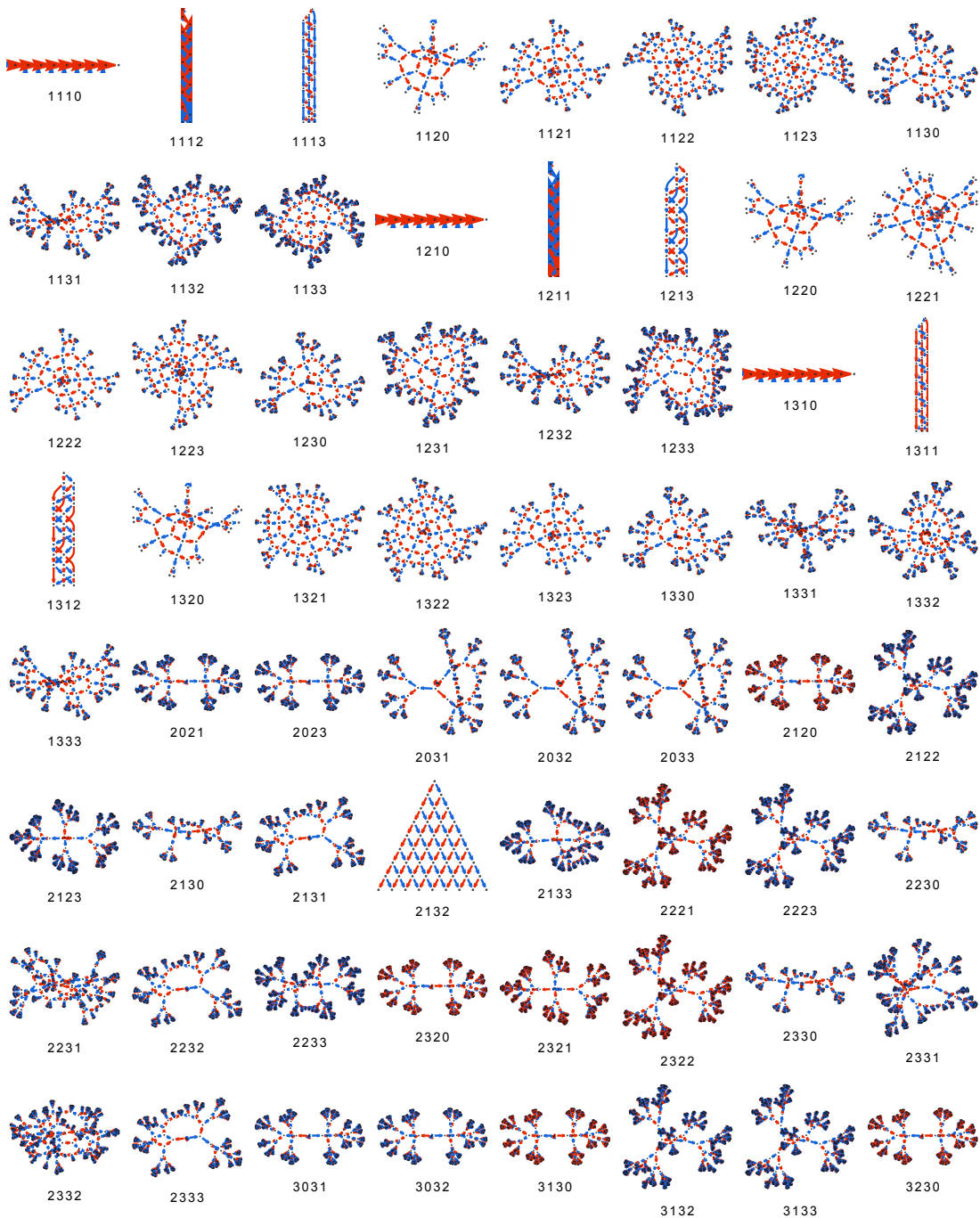

"Ribbons" are seen only when $a = c = 1$. "Simple webs" are seen when $a=1$. "Simple grids" are seen whenever the two cases in the rule commute, i.e.

$(c\,n + d)\,a + b = (a\,n + b)\,c + d$

which occurs whenever

$a\,d - b\,c = \mathrm{db}$

"Simple trees" are seen whenever

$a = c,\ b \neq d$

In other cases there seems to be irregular merging, as in the $n \mapsto \{2\,n + 1, 3\,n + 1\}$ case above. And keeping only nontrivial inequivalent cases these are the results after removing loose ends:

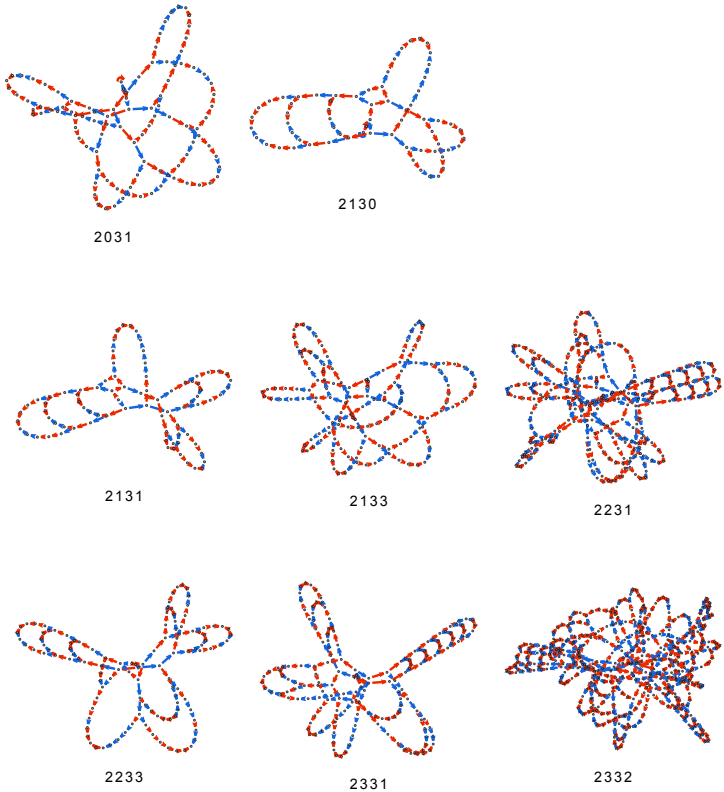

2031

2130

2131

2133

2231

2233

2331

2332

Note that adding another element in the rule can make things significantly more complicated. An example is:

$n \mapsto \{3\,n,\ 2\,n,\ n + 1\}$

After 8 steps this gives

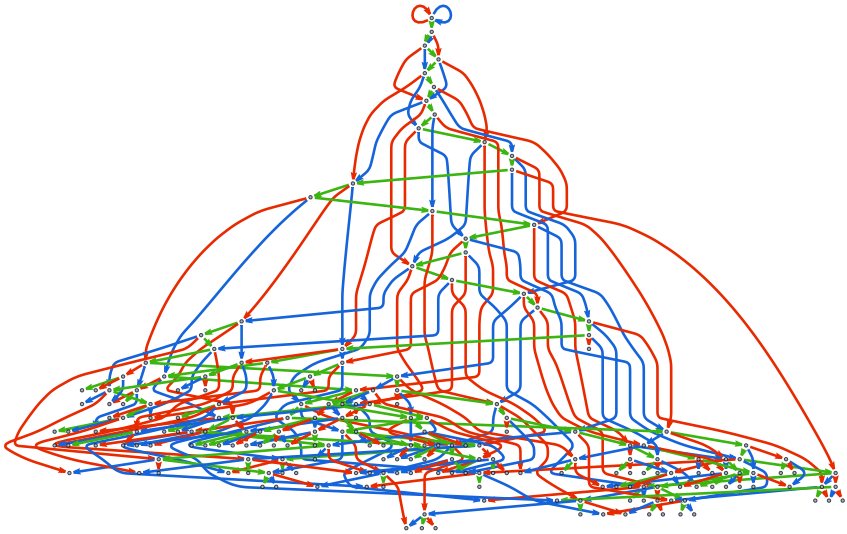

or in another rendering:

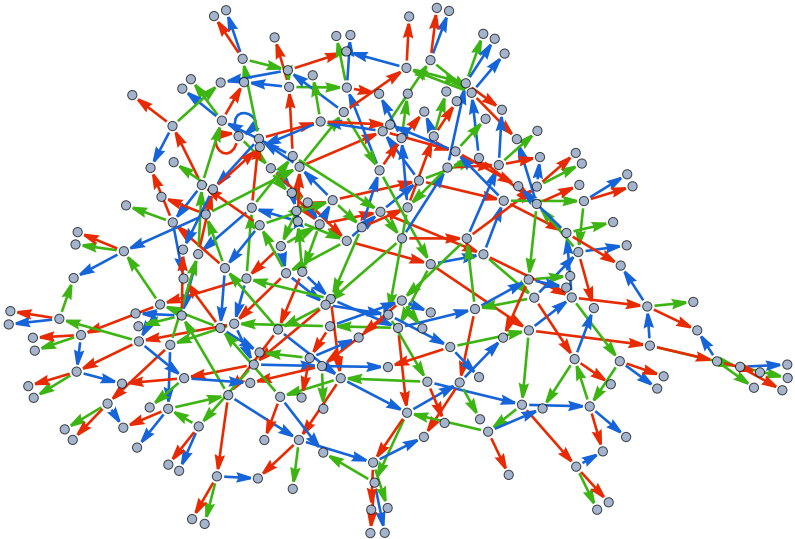

After a few more steps, with "loose ends" removed, one gets the still-rather-unilluminating result (though one that we will discuss further in the next section):

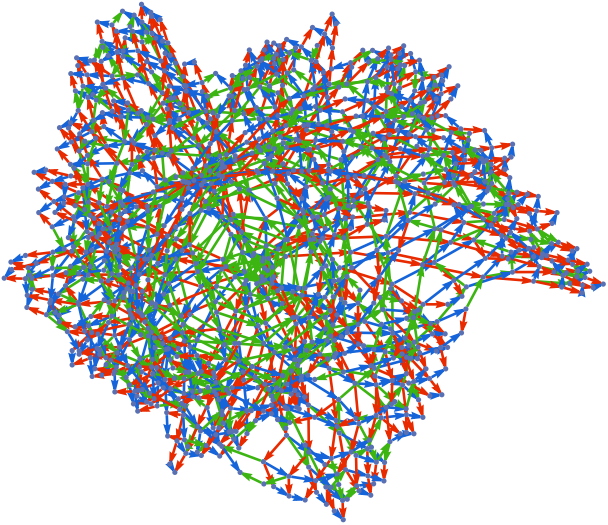

# The Phenomenon of Confluence

Will every branching of paths in the multiway graph eventually merge again? If they do, then the system is confluent (which in this case is equivalent to saying that it's causal invariant—an important property in our Physics Project).

It turns out that all rules of the following forms are confluent:

$n \longmapsto \{a\,n,\, b\,n\}$

$n \longmapsto \{n+a,\, n+b\}$

$n \longmapsto \{a\,n,\, n+b\}$

$n \longmapsto \{a\,n,\, b\,n,\, n+c\}$

But among rules of the form

$n \longmapsto \{a\,n+b,\, c\,n+d\}$

confluence depends on the values of $a$, $b$, $c$ and $d$. When multiway graphs are "simple webs" or "simple grids" there is obvious confluence. And when the graphs are simple trees, there is obviously not confluence.

But what about a case like the rule we discussed above:

$n \longmapsto \{2\,n+1,\, 3\,n+1\}$

We plotted above the "positions" of mergers that occur. But are there "enough" mergers to "rejoin" all branchings?

Here are the first few branchings that occur:

{1 → {3, 4}, 3 → {7, 10}, 4 → {9, 13}, 7 → {15, 22}, 10 → {21, 31}, 9 → {19, 28}}

For the pair 3, 4 one can reach a "merged" end state on the following paths:

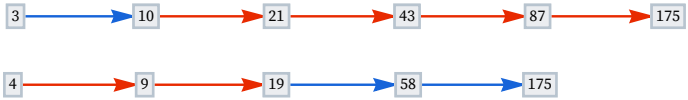

which are embedded in the whole multiway graph (without loose ends) as:

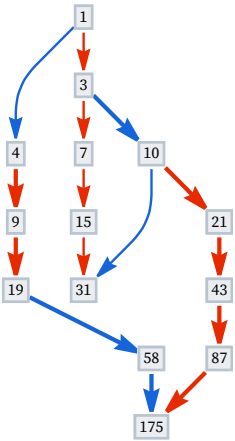

For the pair 9, 13 both eventually reach 177151, but 9 takes 13 steps to do so:

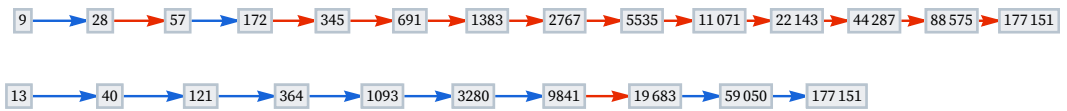

Here's a summary of what we know about what happens with the first few branchings:

| branch pair | merge steps | merge value |
|---|---:|---:|
| {3, 4} | 5 | 175 |
| {7, 10} | 2 | 31 |
| {9, 13} | 13 | 177 151 |
| {15, 22} | 18 | 30 036 991 |
| {21, 31} | >29 | |
| {19, 28} | >29 | |

So what about the total number of branchings and mergings? This is what happens for the first several steps:

| steps | 1 | 2 | 3 | 4 | 5 | 6 | 7 | 8 | 9 | 10 | 11 | 12 | 13 | 14 | 15 | 16 | 17 |
|---|---|---|---|---|---|---|---|---|---|---|---|---|---|---|---|---|---|
| branchings | 1 | 2 | 4 | 8 | 15 | 30 | 59 | 118 | 235 | 468 | 934 | 1866 | 3729 | 7449 | 14887 | 29756 | 59482 |
| mergings | 0 | 0 | 0 | 1 | 0 | 1 | 0 | 1 | 2 | 2 | 2 | 3 | 9 | 11 | 18 | 30 | 41 |

The number of branchings at step $t$ approximates

$$2^{t-1} - 2^{t-5}$$

while the number of mergings seems to grow systematically more slowly, perhaps like $1.5^t$:

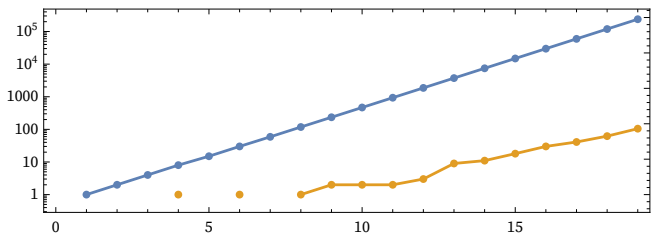

And based on this it seems plausible that the system is not in the end confluent. But how might we show this? And what is the best way to figure out if any particular branch pair (say 21, 31) will ever merge?

One way to look for mergings is just to evolve the multiway graph from each member of the pair, and check if they overlap. But as we can see even for the pair {3, 4} this effectively involves "treeing out" an exponential number of cases:

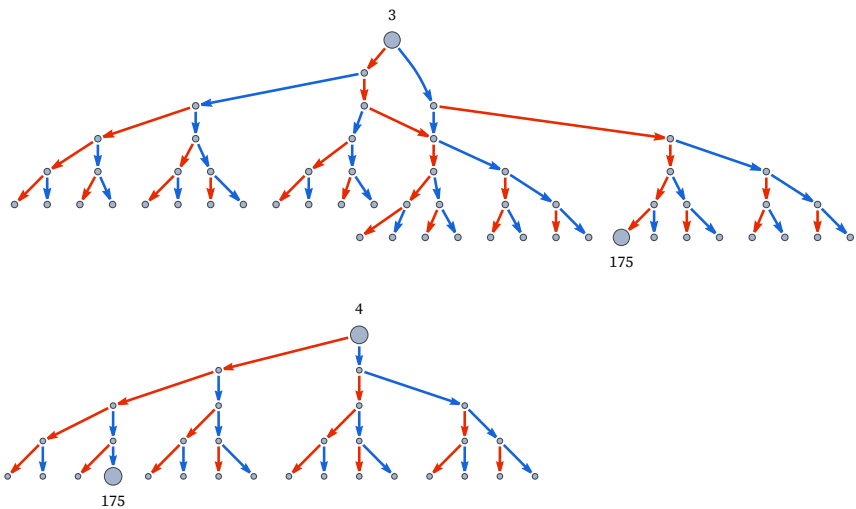

Is there a way to do this more efficiently, or in effect to prune the trees? A notable feature of the original rule is that the numbers it generates always increase at each step. So one thing to do is just to discard all elements at a particular step in one graph that cannot reach the "minimum frontier" in the other graph. But on its own, this leads to only very minor reduction in the size of graph that has to be considered.

To find what is potentially a much more effective "optimization" let's look at some examples of mergings:

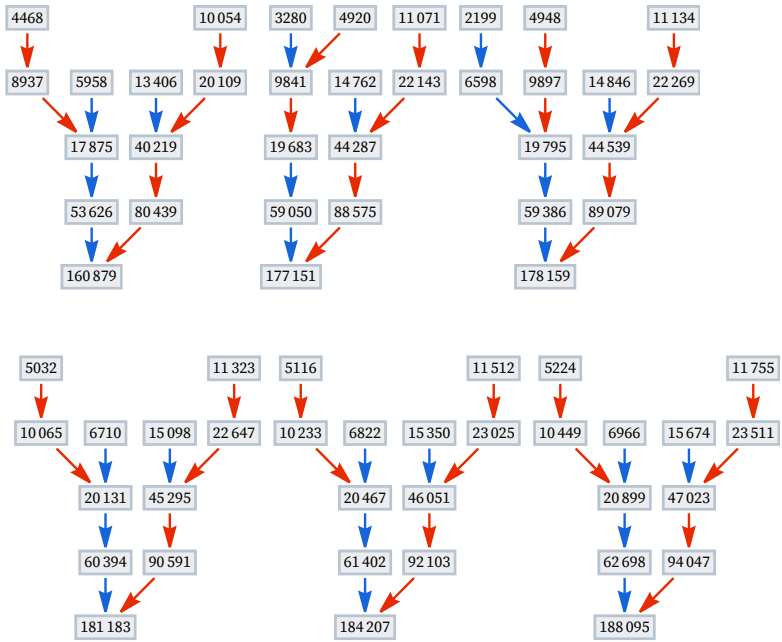

It's clear that the final step has to consist of one application of $n \longmapsto 2\,n+1$ and one of $n \longmapsto 3\,n+1$ (i.e. one red edge and one blue edge). But these examples suggest that there are also further regularities.

At the merging point it must be true that

$$2\,u + 1 = 3\,v + 1$$

for some integers $u$ and $v$. But for this to be true, the merged value (i.e. $2\,u + 1$ or $3\,v + 1$) must for example be equal to 1 mod 2, 3 and 6.

Using the structure one level back we also have:

$$3\,(3\,u + 1) + 1 = 2\,(2\,v + 1) + 1$$

or

$$4 + 9\,u = 3 + 4\,v$$

implying that the merged value must be 3 mod 4, 7 mod 12, 13 mod 18 and 36 mod 31. Additional constraints from going even further back imply in the end that the merged value must have the following pattern of residues:

| modulus | 2 | 3 | 4 | 6 | 8 | 9 | 12 | 16 | 18 | 24 | 36 | 48 | 72 | 144 |
|---------|---|---|---|---|---|---|----|----|----|----|----|----|----|-----|
| residue | 1 | 1 | 3 | 1 | 7 | 4 | 7 | 15 | 13 | 7 | 31 | 31 | 31 | 31 |

But now let's consider the whole system modulo $k$. Then there are just $k$ possible values, and the multiway graph must be finite. For example, for $k = 12$ we get:

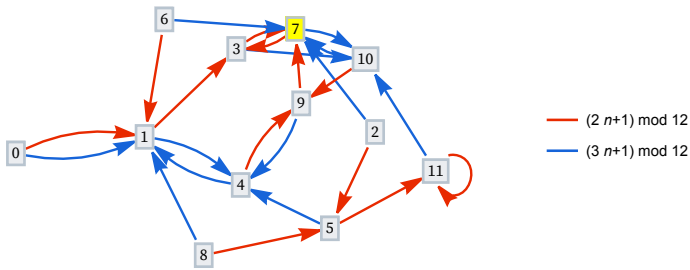

Dropping the "transient parts" leaves just:

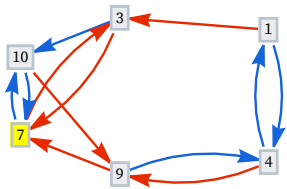

These graphs can be thought of as reductions of the multiway graph (and, conversely, the multiway graph is a covering of them). The graphs can also be thought of as finite automata that define regular languages whose elements are the "2" and "3" transformations that appear on the edges. Any sequence of "2" and "3" transformations that can occur in the multiway graph must then correspond to a valid word in this regular language. But what we have seen is that for certain values of $k$, mergers in the multiway graph always occur at particular ("acceptor") states in the finite automata.

In the case $k = 12$, every merger occurs at the 7 state. But by tracing possible paths in the finite automaton we now can read off what sequences of transformations can lead to a merger:

{22, 33}
{222, 233, 322, 333}
{2222, 2233, 2322, 3233, 3322, 3333}
{22222, 22233, 22322, 23322, 23333, 32222, 32233, 32322, 33233, 33322, 33333}

And what's notable is that only a certain fraction of all $2^m$ possible sequences of length $m$ can occur; asymptotically, about 28%.

The most stringent analogous constraints come from the *k*=144 graph:

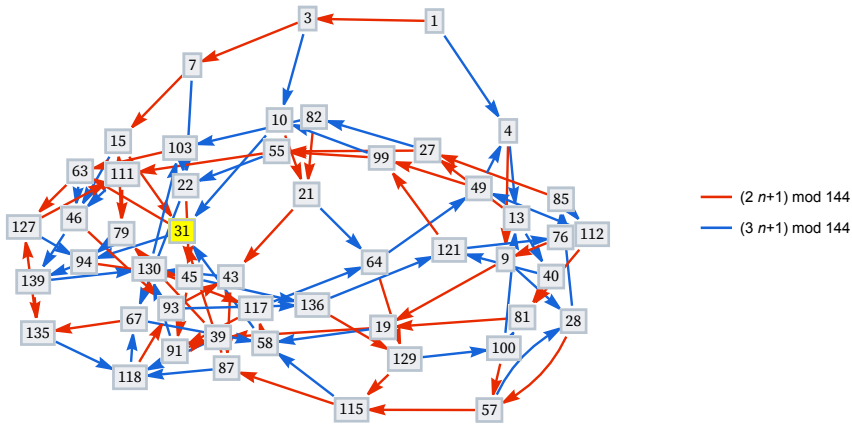

And we see that even for sequences of length 3 fewer are allowed than from the *k* = 12 graph:

{22, 33}
{222, 233, 333}
{2222, 2233, 3233, 3333}
{22222, 23333, 32222, 32233, 33233, 33333}

Asymptotically the number of allowed sequences is about 3% of the $2^m$ possible. And so the conclusion is that if one wants to find mergings in the multiway graph it's not necessary to tree out all possible sequences of transformations; one only needs at most the 30× smaller number of sequences "accepted by the mod-144 finite automaton". It's possible to do a little better than this, by looking not just at sequences allowed by the finite automaton for a particular *k*, but at finite automata for a collection of values of *k* (say as in the table above).

But while these techniques deliver significant practical speedups they do not seem to significantly alter the asymptotic resources needed. So what will it take to determine whether the pair {21, 31} ever merges?

I don't know. And for example I don't know any way to find an upper bound on the number of steps after which we'd be able to say "if it hasn't merged yet, it never will". I'm sure that if we look at different branch pairs, there will be tricks for particular cases. But I suspect that the general problem of determining merging will show computational irreducibility, and that for example there will be no fundamentally better way to determine whether a particular branch pair has merged after *t* steps than by essentially enumerating every possible evolution for that number of steps.

But if this is the case, it means that the general infinite-time question of whether a branch pair will merge is undecidable—and can never be guaranteed to be answerable with a bounded amount of computational effort. It's a lower bar to ask whether the question can be answered using a finite proof in, say, Peano arithmetic. And I think it's very likely that the overall question of whether all branch pairs merge—so that the system is confluent—is a statement that can never, for example, be established purely within Peano arithmetic. There are quite a few other candidates for the "simplest 'numerical' statement independent of Peano arithmetic". But it seems at least conceivable that this one might be more accessible to proof than most.

It's worth mentioning, by the way, that (as we have seen extensively in the Physics Project) the presence of confluence does not imply that a multiway system must show simple overall behavior. Consider for example the rule (also discussed at the end of the previous section):

$n \longmapsto \{3\,n,\, 2\,n,\, n + 1\}$

Running for a few more steps, removing loose ends and rendering in 3D gives:

But despite this complexity, this is a confluent rule. It's already an indication of this that mergings pretty much "keep up" with branchings in this multiway system:

| steps | 1 | 2 | 3 | 4 | 5 | 6 | 7 | 8 | 9 | 10 | 11 | 12 | 13 | 14 |
|---|---|---|---|---|---|---|---|---|---|---|---|---|---|---|
| branchings | 2 | 3 | 7 | 14 | 27 | 53 | 106 | 210 | 422 | 840 | 1686 | 3364 | 6733 | 13469 |
| mergings | 3 | 2 | 5 | 13 | 24 | 47 | 92 | 180 | 366 | 722 | 1456 | 2903 | 5812 | 11596 |

The first few branchings (now all 3-way) are:

$\{1 \rightarrow \{2, 2, 3\}, 2 \rightarrow \{3, 4, 6\}, 3 \rightarrow \{4, 6, 9\}, 4 \rightarrow \{5, 8, 12\}, 6 \rightarrow \{7, 12, 18\}, 9 \rightarrow \{10, 18, 27\}\}$

All the pairs here merge (often somewhat degenerately) in just a few steps. Here are examples of how they work:

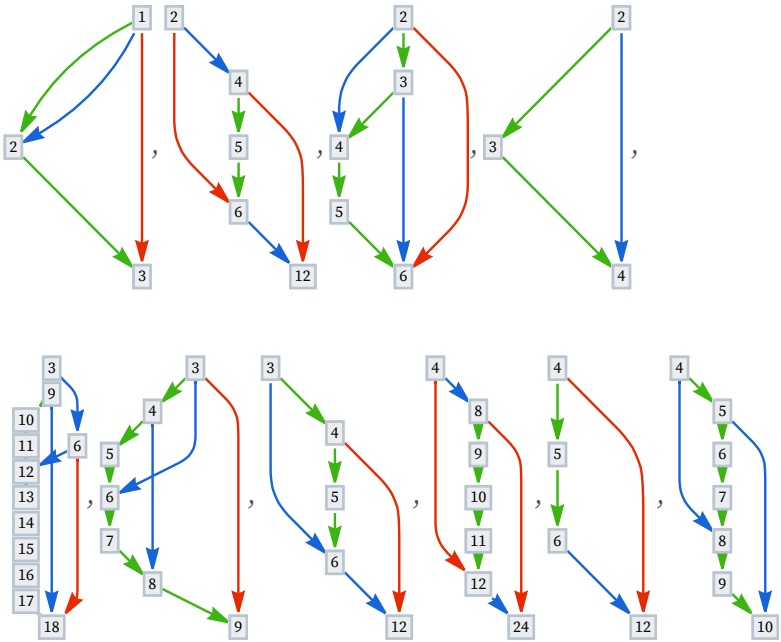

# Branchial Space and Numerical Value Space

Consider the first few steps of the rule

$n \longmapsto \{3\,n, 2\,n, n + 1\}$

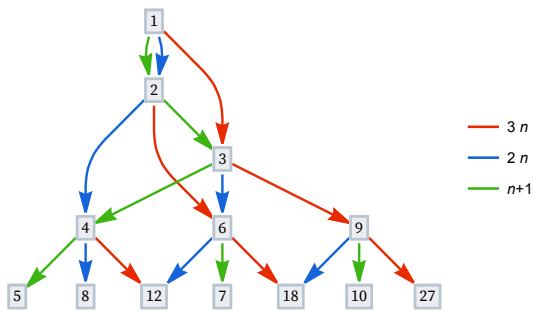

At each "layer" we can form a branchial graph by joining nodes that have common ancestors on the step before:

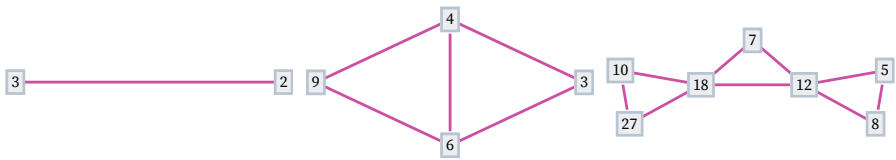

Continuing for a few more steps we get:

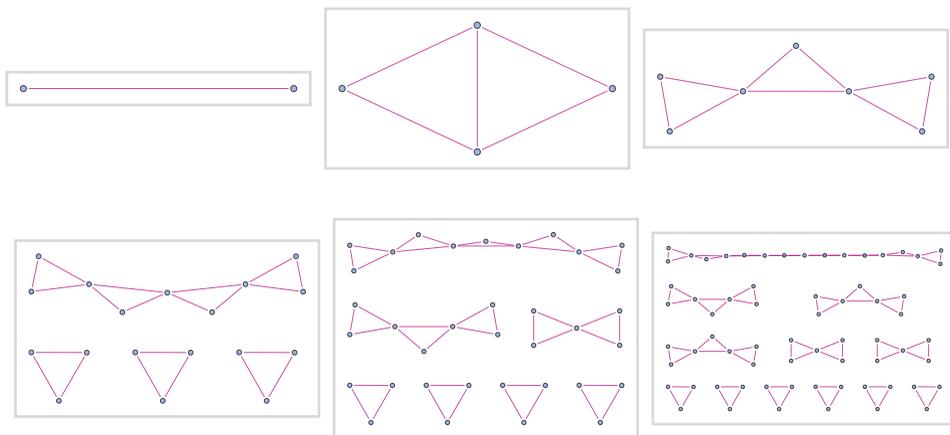

We can imagine (as we do in our Physics Project) that in an appropriate (if rather subtle) limit such branchial graphs can be thought of as defining a "branchial space" in which each node has a definite position. (One of many subtleties is that the particular branchial graphs we show here are specific to the particular "layering" of the multiway graph that we've used; different foliations would give different results.)

But whereas in our Physics Project and many other applications of the multicomputational paradigm the only real way to define "positions" for nodes in the multiway graph is through something like branchial space, there is a much more direct approach that can be taken in multiway systems based on numbers—because every node is labeled by a number which one can imagine directly using as a coordinate.

As an example, let's take the multiway graph above, and make the horizontal position of each node be determined by its value:

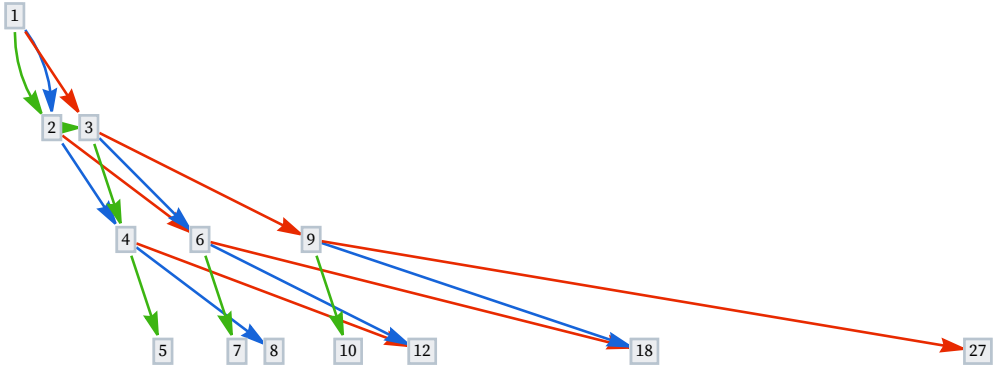

Or, better, by the log of its value:

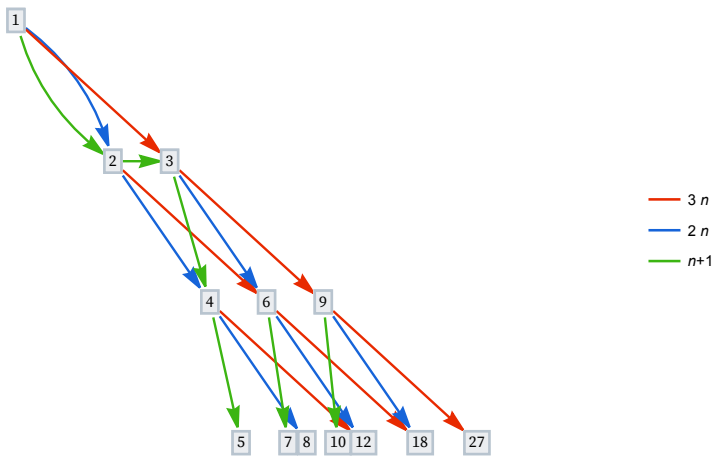

Continuing for more steps, we get:

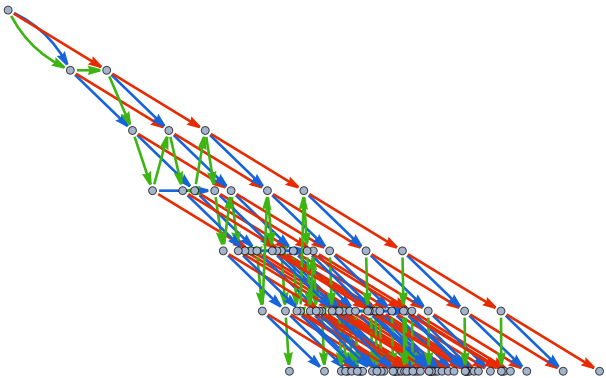

Now, for example, we can ask—given the particular choice of layers we have made here— what the distribution of (logarithmic) values reached on successive layers will be, and one finds that the results converge quite quickly:

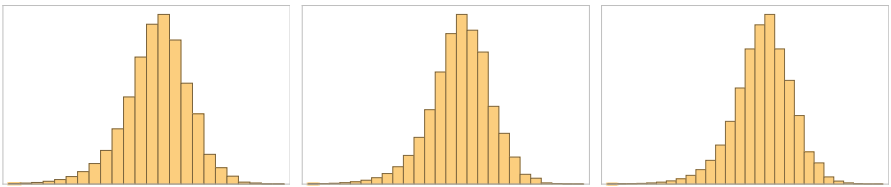

(By the way, in these results we've not included "path weights", which determine how many different paths lead from the initial number to a particular result. In the example shown, including path weights doesn't make a difference to the form of the final result.)

So what is the correspondence between the layout of nodes in "branchial space" and in "numerical value space"? Here's what happens if we lay out a branchial graph using (logarithmic) numerical value as $x$ coordinate:

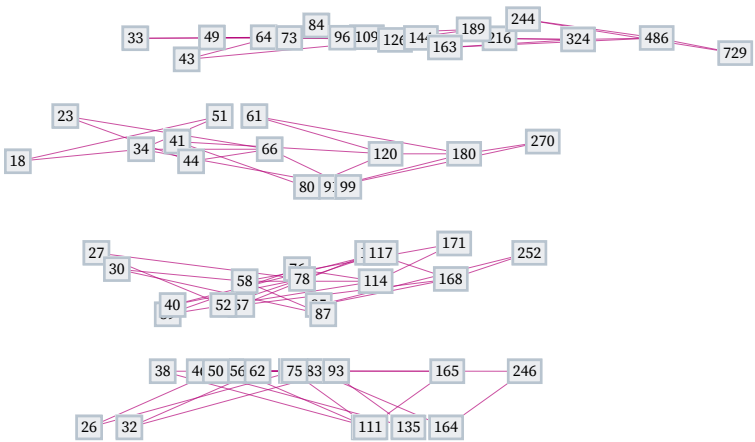

Perhaps more useful is to plot branchial distance versus (logarithmic) numerical distance for every pair of connected nodes at a particular layer:

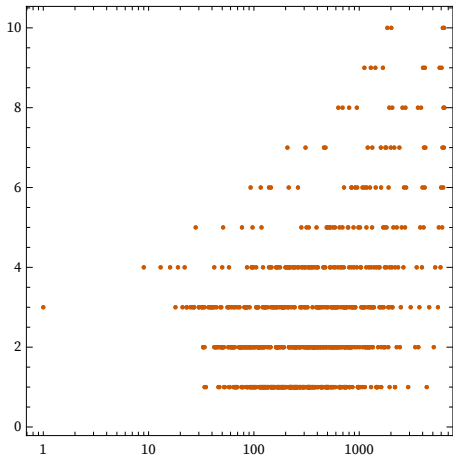

And at least in this case, there is perhaps a slight correlation to be seen.

## Negative Numbers

The rules we've considered so far all involve only non-negative numbers. What happens if we include negative numbers? Generally the results are very similar to those with non-negative numbers. For example:

$n \longmapsto \{2\,n,\ n-1\}$

just gives

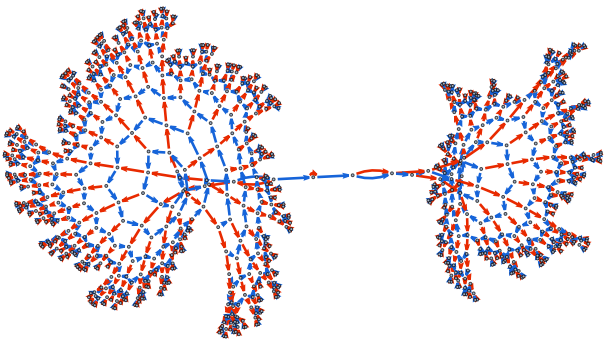

in which there is effectively both a "positive" and "negative" "web".

A rule like

$n \longmapsto \{2\,n+1,\ 3\,n-1\}$

turns out to yield essentially only positive numbers, yielding after removing loose ends

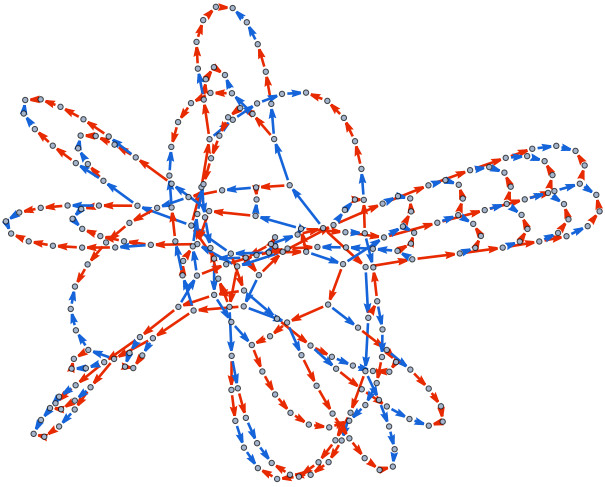

$n \longmapsto \{2\,n + 1, 1 - 3\,n\}$

gives a more balanced collection of positive and negative numbers (with positive numbers indicated by dark nodes), but the final graph is still quite similar:

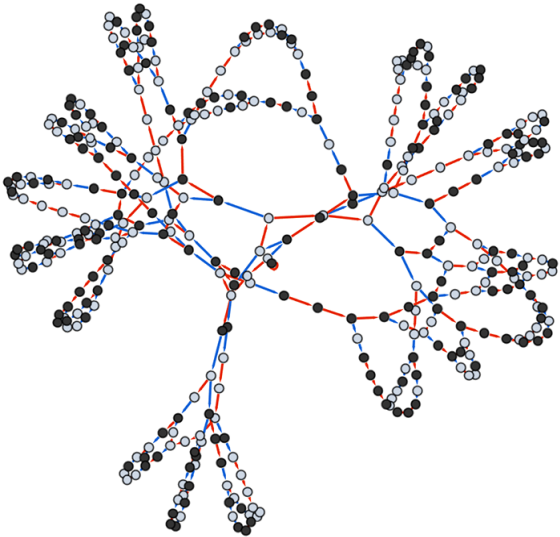

## "Floor" and Related Rules

So far we've considered only rules based on ordinary arithmetic functions. As a first example of going beyond that, consider the rule:

$n \longmapsto \{n + 1, \mathsf{Floor}[\dfrac{n}{2}]\}$

Running this for 50 steps we get:

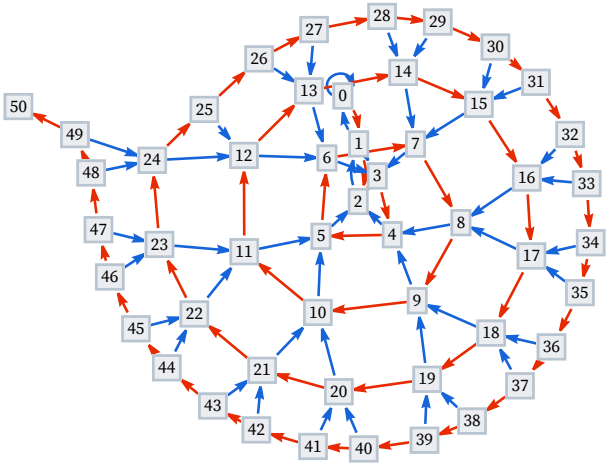

A notable feature here is that only one "fresh" node is added at each step—and the whole thing grows like a Fermat spiral. After 250 steps the multiway graph has the form

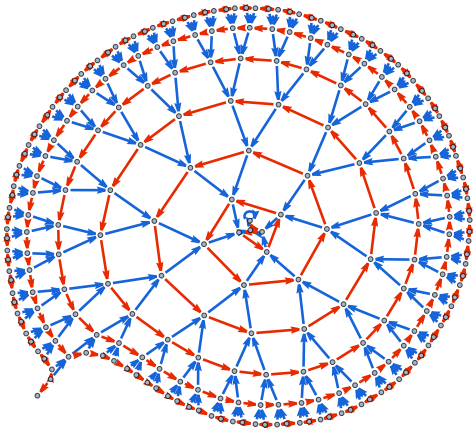

which we can readily see is essentially a "binary tree superimposed on a spiral".

Dividing by 3 instead of 2 makes it a ternary tree:

$$n \mapsto \{n + 1, \ \mathsf{Floor}[\tfrac{n}{3}]\}$$

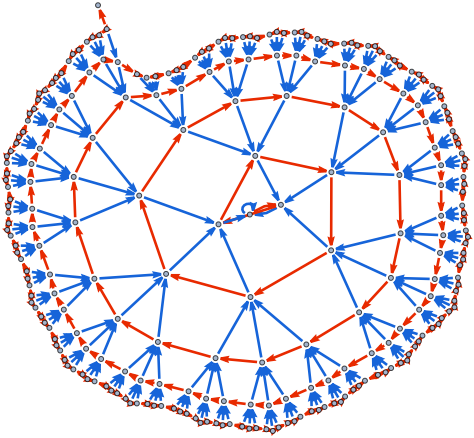

Using **Round** instead of **Floor** gives a mixed binary and ternary tree:

$$n \mapsto \{n+1, \text{Round}[\frac{n}{2}]\}$$

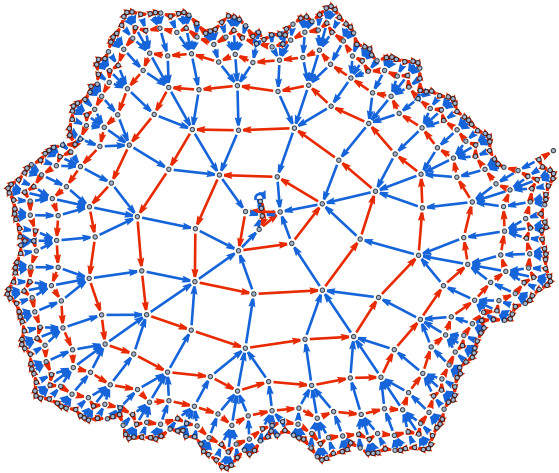

What about rules of the form:

$$n \mapsto \{a\,n, \text{Floor}[\frac{n}{2}]\}$$

Here are the results for a few values of *a*:

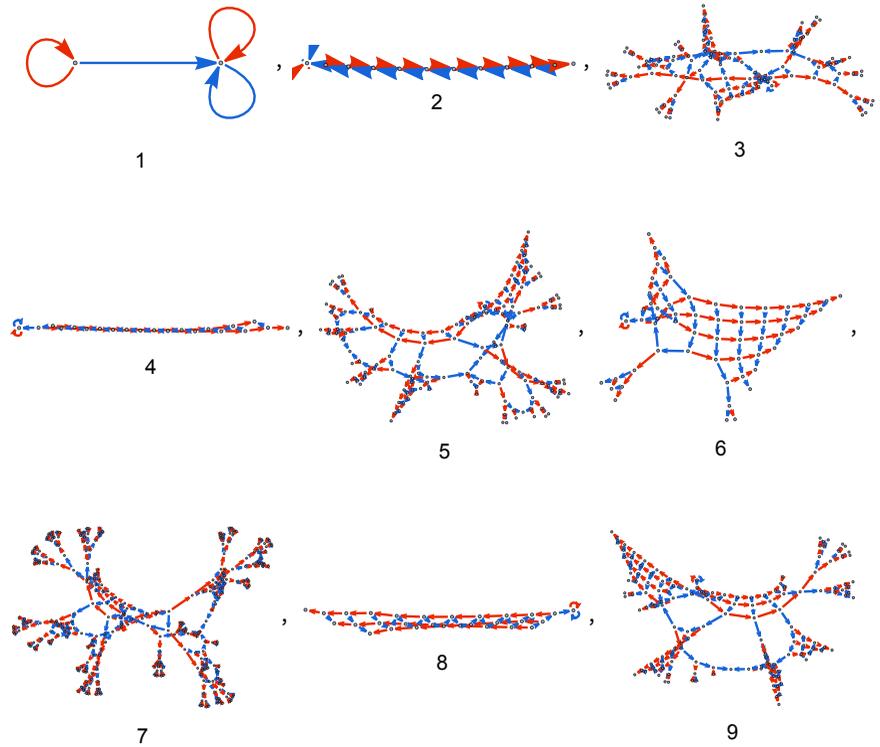

Continuing

$$n \mapsto \{3\,n,\; \mathrm{Floor}[\frac{n}{2}]\}$$

for more steps we get:

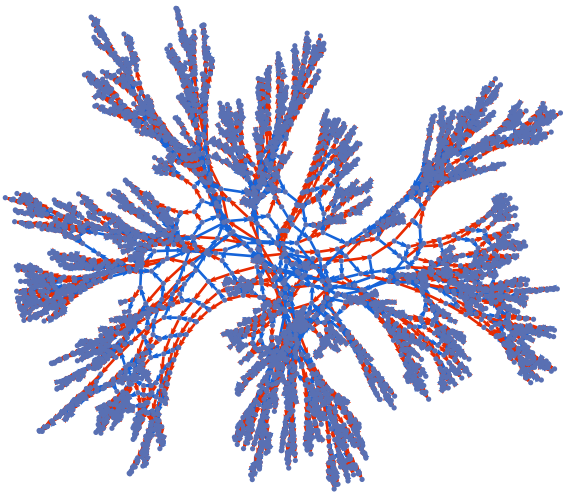

$$n \mapsto \{6\,n,\ \text{Floor}[\frac{n}{2}]\}$$

has far fewer "loose ends":

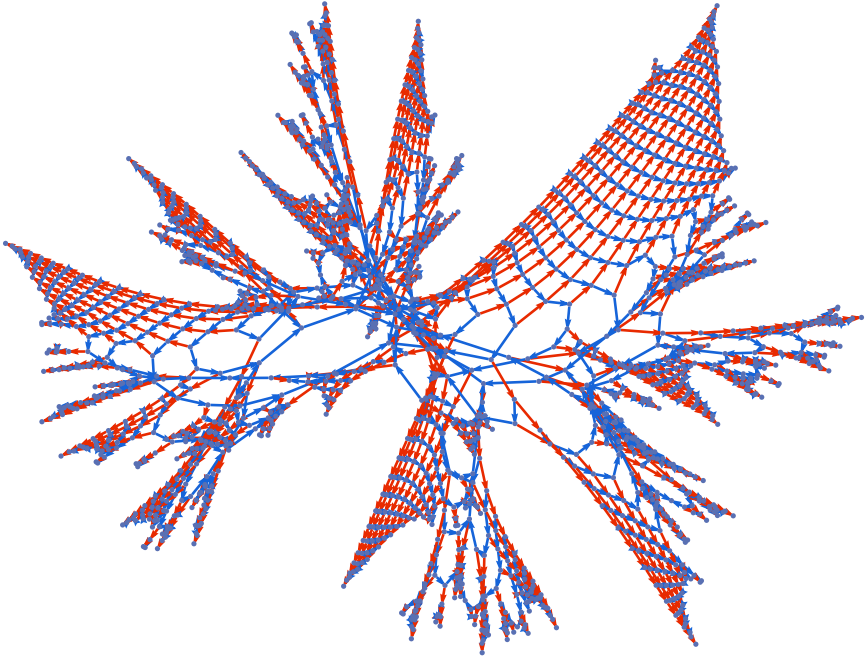

What are the "grid patches"? Picking out some of the patches we can see they're places where a number that can be "halved a lot" appears—and just like in our pure multiplication rules above, $\frac{n}{2}$ and $3n$ represent commuting operations that form a grid:

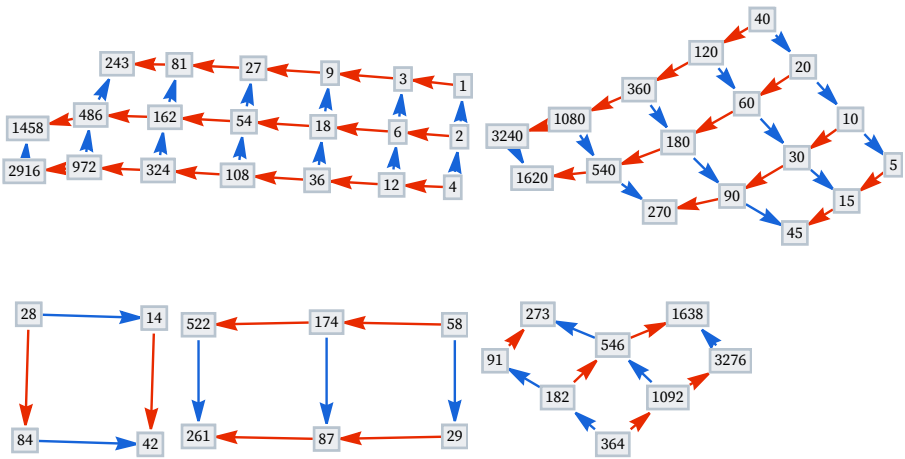

# Conditional Division, Inverse Iterations and the 3n+1 Problem

Including **Floor**[$n/2$] is a bit like having different functions for even and odd $n$. What happens if we do this more explicitly? Consider for example

$n \longmapsto$ If[EvenQ[$n$], {$n\,/\,2$, $n + 1$}, {$n + 1$}]

The result is essentially identical to the **Floor** case:

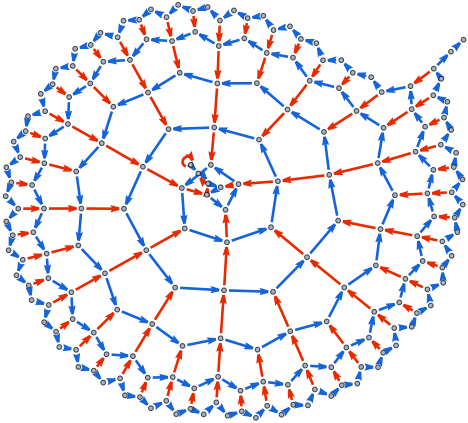

Here are a couple of other cases, at least qualitatively similar to what we've seen before:

$n \longmapsto$ If[EvenQ[$n$], {$n\,/\,2$}, {$3\,n$, $n + 1$}]

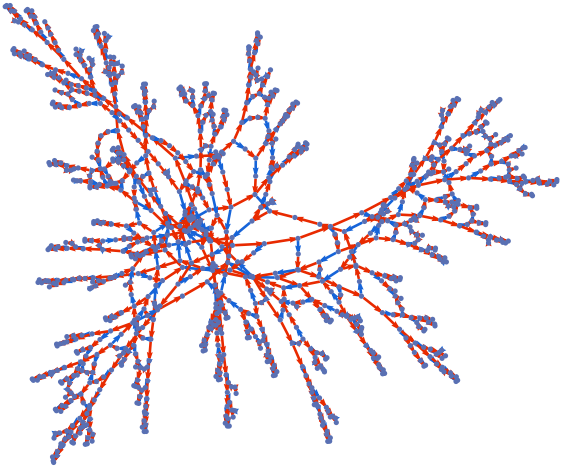

$n \mapsto \text{If}[\text{EvenQ}[n], \{n \, / \, 2, \, 2\,n\}, \{n+2, \, 2\,n\}]$

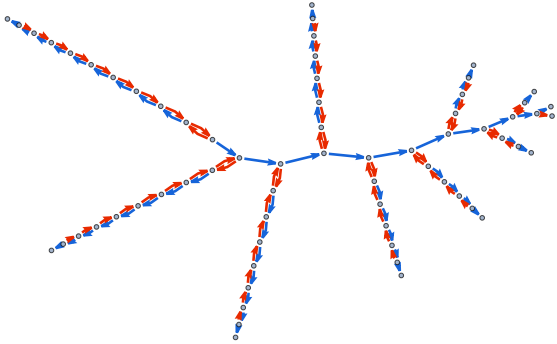

But now consider as we did at the beginning:

$n \mapsto \{2\,n, \, n+1\}$

What is the inverse of this? One can think of it as being

$n \mapsto \text{If}[\text{EvenQ}[n], \{n \, / \, 2, \, n-1\}, \{n-1\}]$

or

$n \mapsto \text{Select}[\{n \, / \, 2, \, n-1\}, \text{IntegerQ}]$

which gives for example

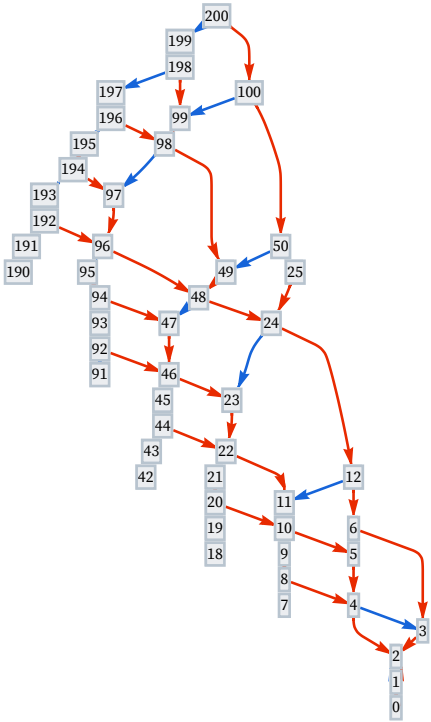

or continuing for longer:

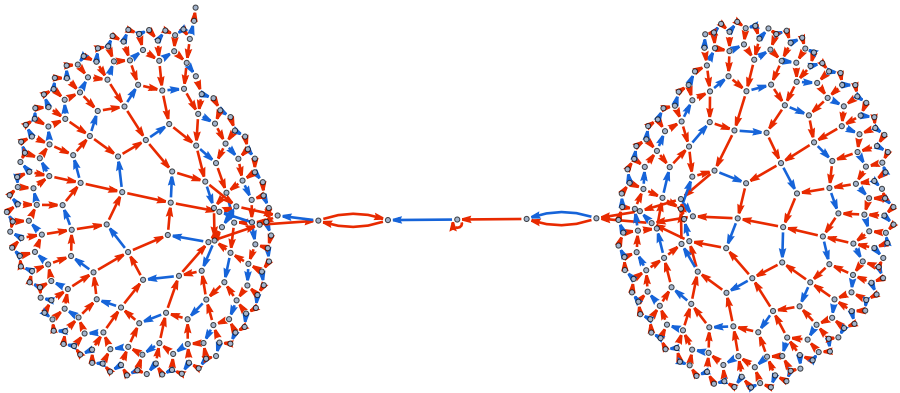

How about

$n \longmapsto \{2\,n + 1,\, 3\,n + 1\}$

Now the "inverse" is:

$n \longmapsto \mathsf{Select}[\{(n-1)\,/\,2,\,(n-1)\,/\,3\},\,\mathsf{IntegerQ}]$

But in this case since most numbers are not reached in the original iteration, most "don't have inverses". However, picking an initial number like 4495, which happens to be a merge point, yields:

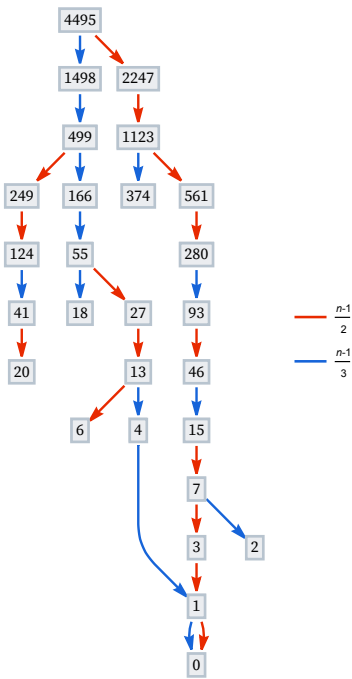

Note that this "inverse iteration" always monotonically decreases towards 0—reaching it in at most $\log_2(n)$ steps.

But now we can compare with the well-known *3n + 1* problem, defined by the "singleway" iteration:

$n \mapsto \mathsf{If[EvenQ}[n], n / 2, 3 n + 1]$

And while in this case the intermediate numbers sometimes increase, all known initial conditions eventually evolve to a simple cycle:

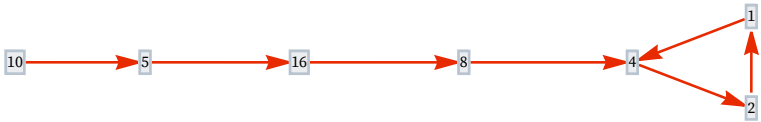

But now we can "invert" the problem, by considering the rule:

$n \mapsto \mathsf{Select}[\{2\,n,\ (n-1)/3\},\ \mathsf{IntegerQ}]$

equivalent to

$n \mapsto \mathsf{If[Mod}[n, 3] = 1,\ \{2\,n,\ (n-1)/3\},\ \{2\,n\}]$

which gives after 10 steps:

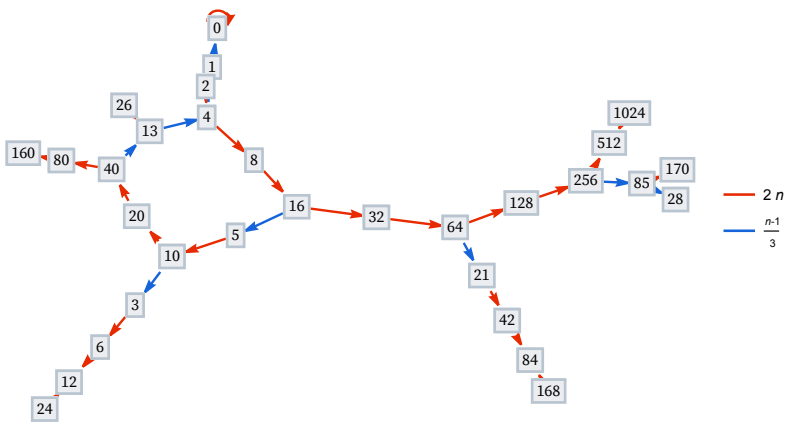

Continuing this to 25 steps one gets:

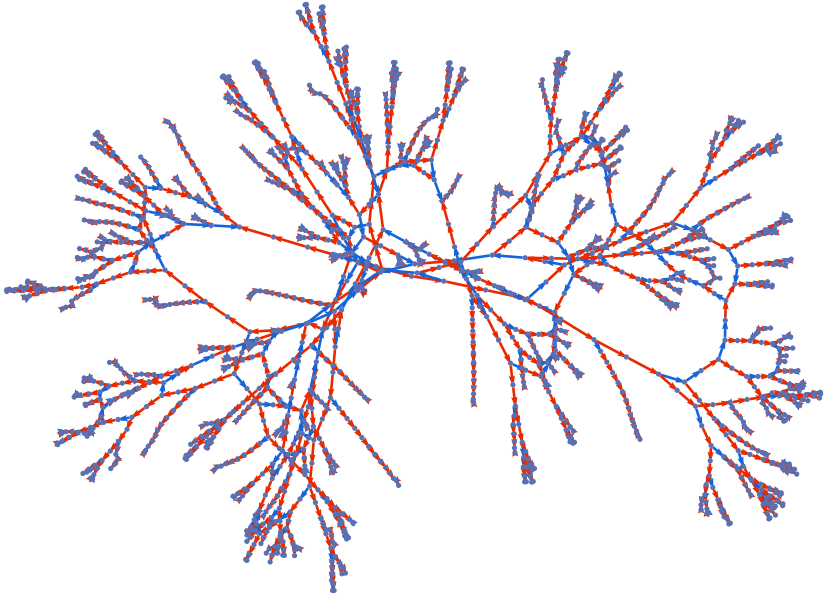

Removing loose ends this then becomes:

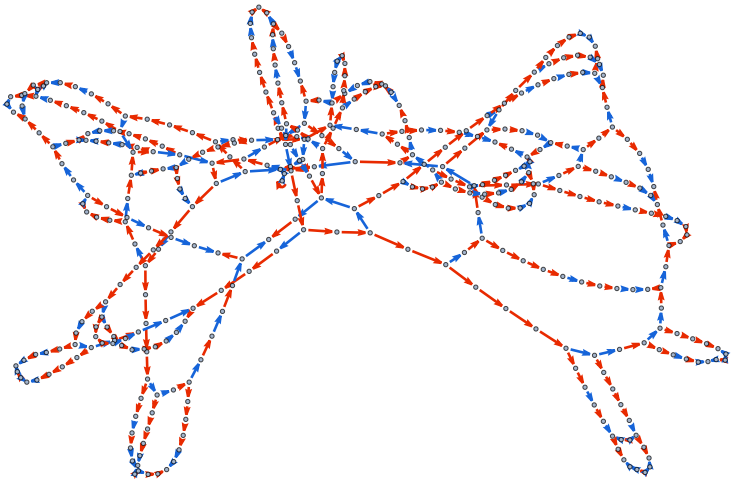

or after more steps, and rendered in 3D:

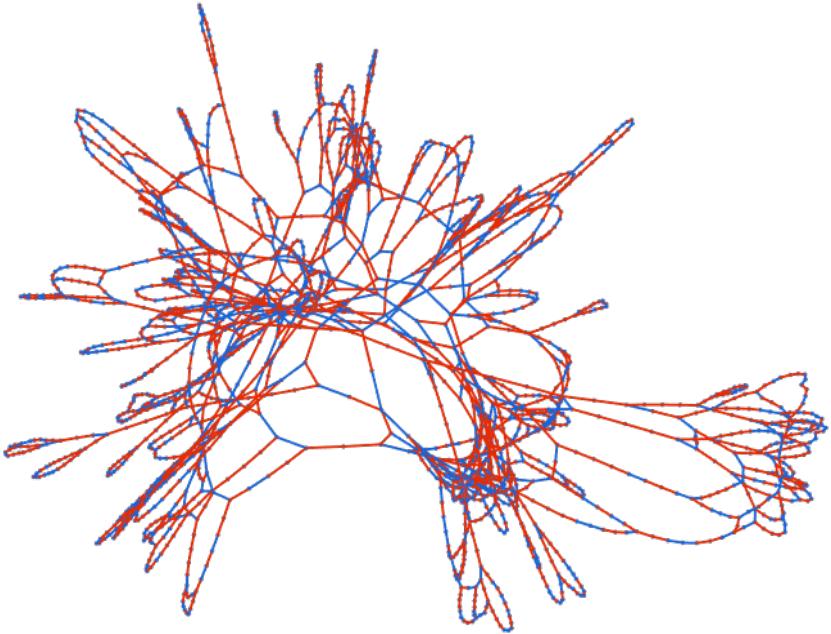

The *3n + 1* problem now asks whether as the multiway graph is built, it will eventually include every number. But from a multicomputational point of view there are new questions to ask—like whether the "inverse-*3n + 1*-problem" multiway system is confluent.

The first few branchings in the multiway graph in this case are

$\{1 \to \{0, 2\}, 4 \to \{1, 8\}, 7 \to \{2, 14\}, 10 \to \{3, 20\}, 13 \to \{4, 26\},$
$\quad 16 \to \{5, 32\}, 34 \to \{11, 68\}, 37 \to \{12, 74\}, 40 \to \{13, 80\}, 46 \to \{15, 92\}\}$

and all of these re-merge after at most 13 steps. The total number of branchings and mergings on successive steps is given by:

| steps | 1 | 2 | 3 | 4 | 5 | 6 | 7 | 8 | 9 | 10 | 11 | 12 | 13 | 14 | 15 | 16 | 17 | 18 | 19 | 20 | 21 | 22 | 23 | 24 | 25 | 26 | 27 | 28 | 29 | 30 |
|---|---|---|---|---|---|---|---|---|---|---|---|---|---|---|---|---|---|---|---|---|---|---|---|---|---|---|---|---|---|---|
| branchings | 0 | 1 | 0 | 1 | 0 | 2 | 0 | 2 | 2 | 3 | 2 | 5 | 5 | 7 | 9 | 14 | 14 | 18 | 27 | 34 | 43 | 56 | 73 | 93 | 118 | 159 | 201 | 260 | 335 | 437 |
| mergings | 1 | 0 | 0 | 0 | 0 | 0 | 0 | 0 | 1 | 0 | 0 | 0 | 1 | 0 | 0 | 5 | 1 | 1 | 2 | 6 | 4 | 5 | 13 | 10 | 13 | 17 | 29 | 31 | 38 | 54 |

Including more steps one gets

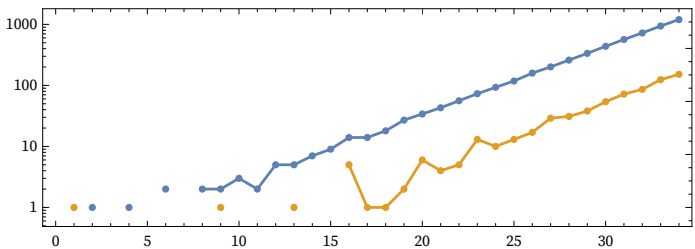

which suggests that there is indeed confluence in this case—though, like for the problem of termination in the original *3n + 1* problem, it may be extremely difficult to determine this for sure.

## Other Kinds of Rules

All the rules we've used so far are—up to conditionals—fundamentally "linear". But we can also consider "polynomial" rules. With pure powers, as in

$n \longmapsto \{n^2, n^3\}$

the multiway graph is just the one associated with the addition of exponents:

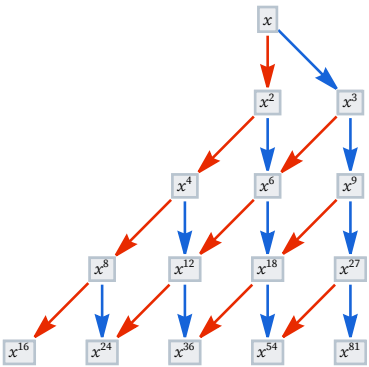

In a case like

$n \longmapsto \{n^2 + 1, 2\,n\}$

the graph is a pure tree

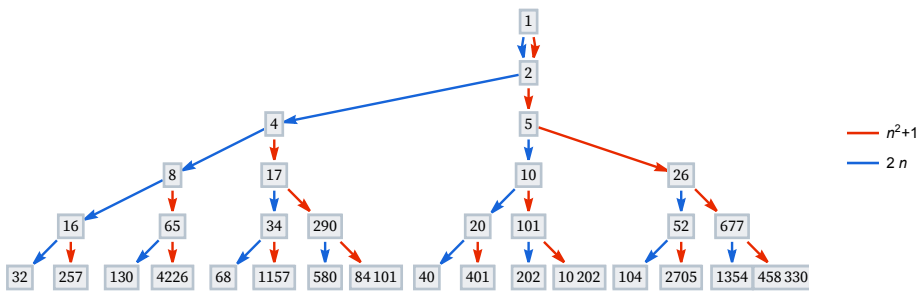

while in a case like

$n \longmapsto \{n^2 + 1, n + 1\}$

there is "early merging", followed by a pure tree:

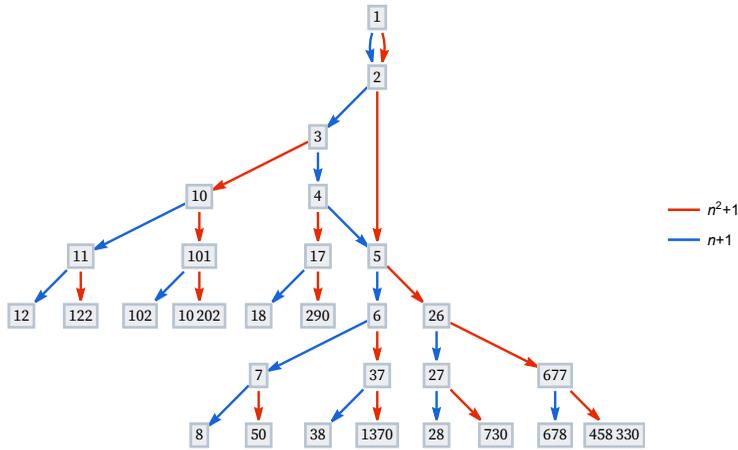

There are also cases like

$n \mapsto \{n^2 + n + 1, \, n + 1\}$

which lead to "continued merging"

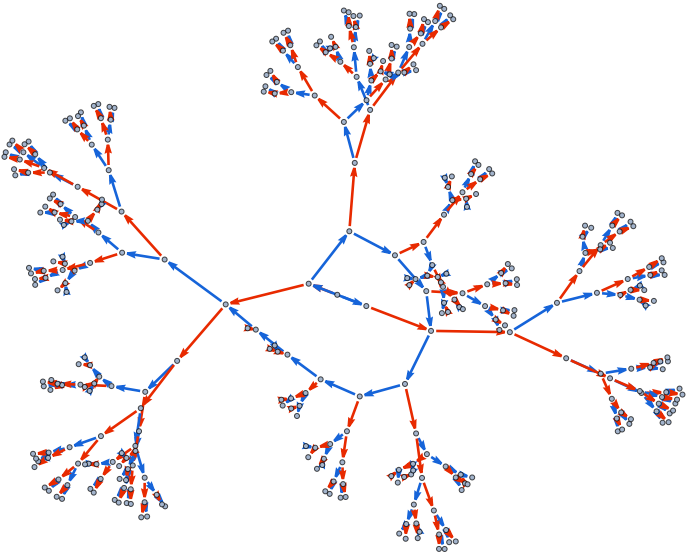

but when loose ends are removed, they are revealed to behave in rather simple ways:

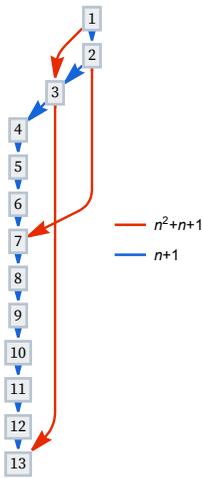

In a case like

$$n \mapsto \{(n + 1)^2, \, 2\,n\}$$

however, there is at least slightly more complicated merging (shown here after removing loose ends):

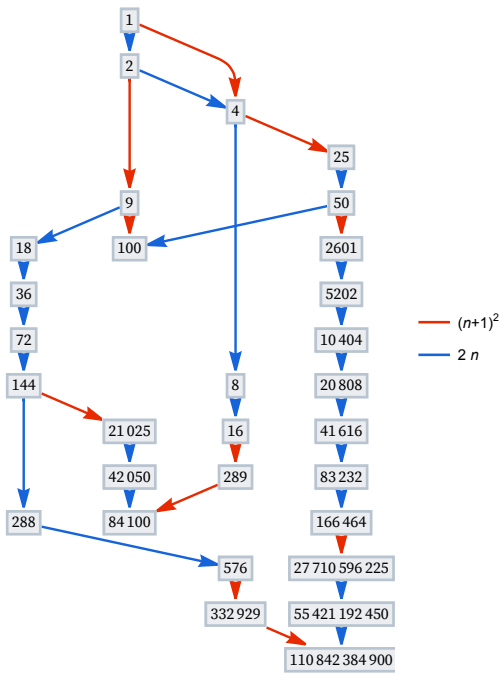

If we include negative numbers we find cases like:

$n \mapsto \{n^2 - n - 1, \, n - 1\}$

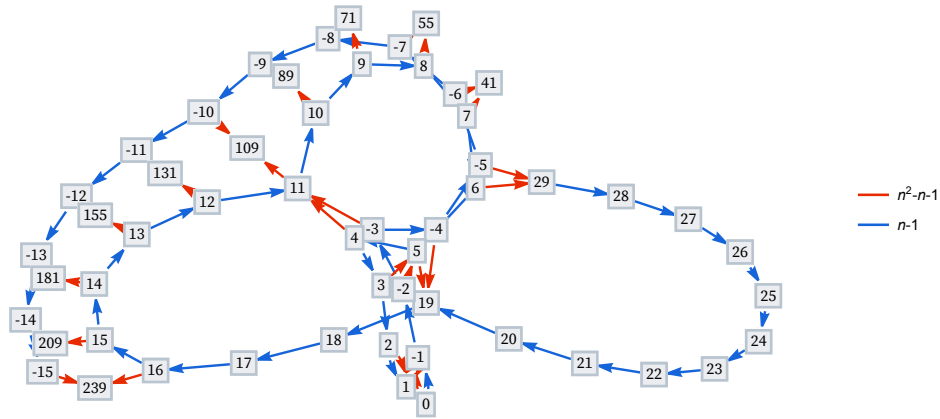

But in other "polynomial" cases one tends to get only trees; a merging corresponds to a solution to a high-degree Diophantine equation, and things like the ABC conjecture tend to suggest that very few of these exist.

Returning to the "linear" case, we can consider—as we did above—multiway graphs mod $k$. Such graphs always have just $k$ nodes. And in a case like

$n \mapsto \text{Mod}[\{n + 1, 10\,n\}, 7]$

with graph

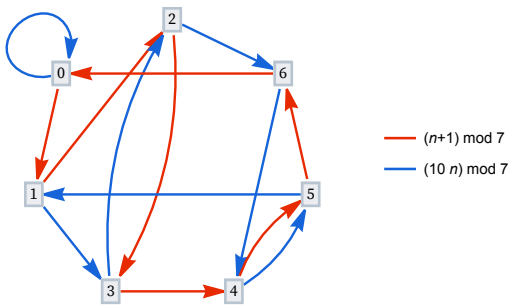

they have a simple interpretation—as "remainder graphs" which one can use to compute a given input number $n$ mod $k$. Consider for example the number 867, with digits 8, 6 and 7. Start at the 0 node. Follow 8 red arrows, followed by a blue one, thus reaching node 3. Then follow 6 red arrows, followed by blue. Then 7 red arrows, followed by blue. The node that one ends up on by this procedure is exactly the remainder. And in this case it is node 6, indicating that $\text{Mod}[867, 7]$ is 6.

Not too surprisingly, there is a definite structure to such remainder graphs. Here is the sequence of "binary remainder graphs" generated from the rule

$n \mapsto \mathsf{Mod}[\{n + 1, 2\,n\}, k]$

for successive values of *k*:

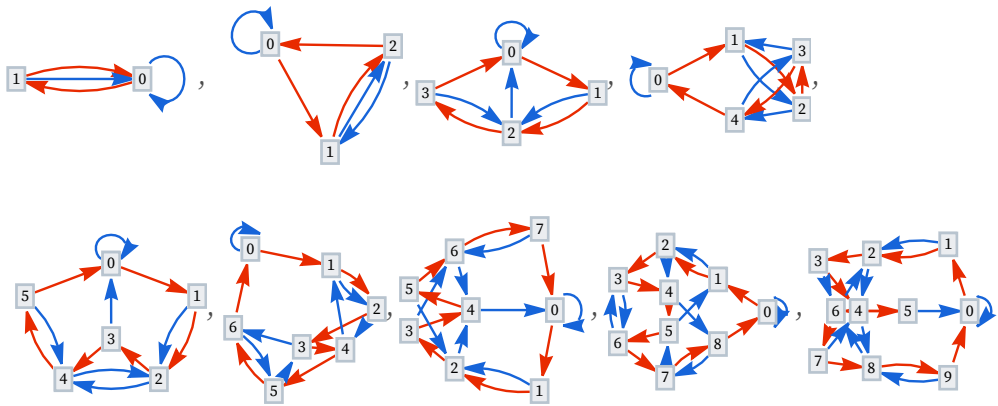

Continuing a number-theoretical theme, we may note that the familiar "divisor graph" for a number can be considered as a multiway graph generated by the rule:

$n \mapsto \mathsf{Most}[\mathsf{Divisors}[n]]$

Here's an example for 100:

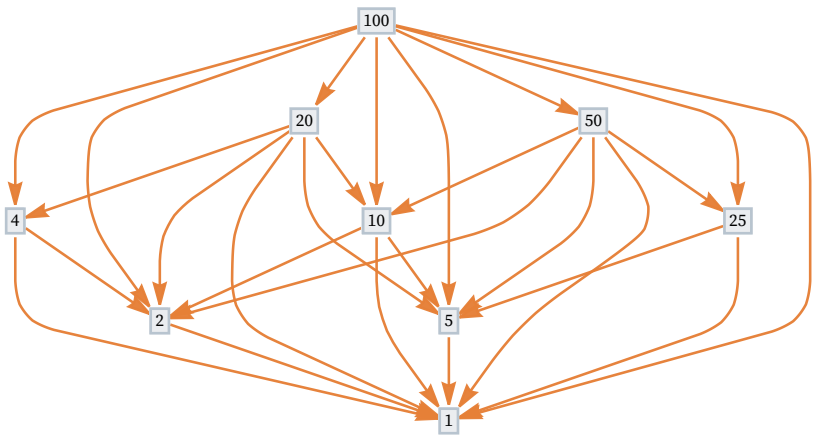

Transitive reduction gives a graph which in this case is essentially a grid:

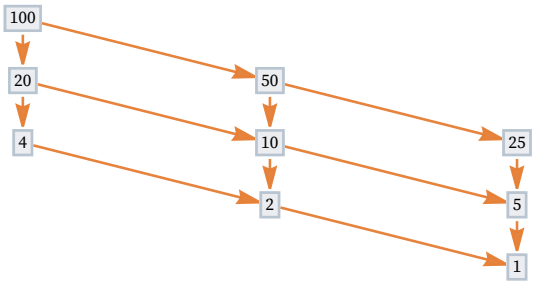

Other initial numbers can give more complicated graphs

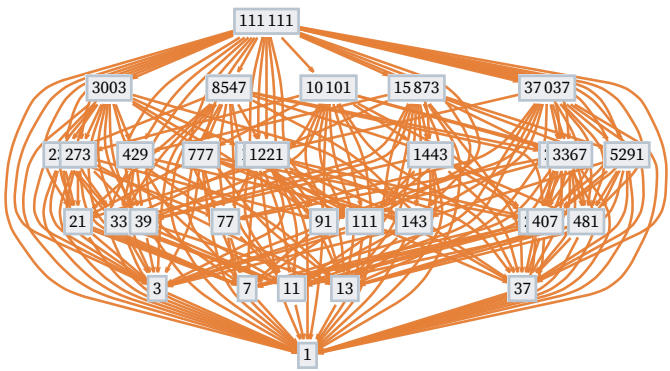

but in general the transitive reduction is essentially a grid graph of dimension `PrimeNu[`*n*`]`:

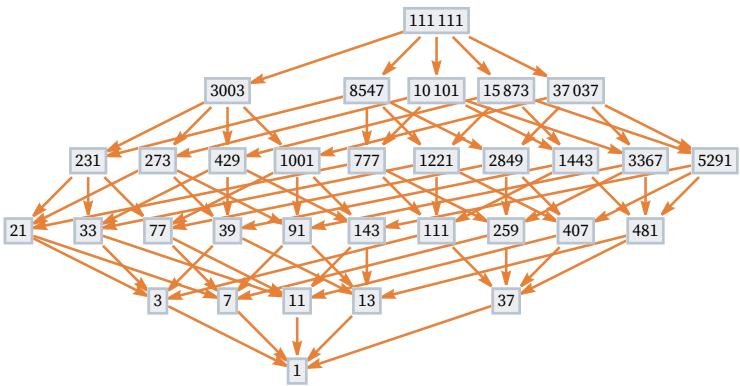

As an alternative to looking at divisors, we can look, for example, at a rule which transforms any number to the list of numbers relatively prime to it:

$n \longmapsto$ Select[Range[2, $n$ – 1], CoprimeQ[♯, $n$] &]

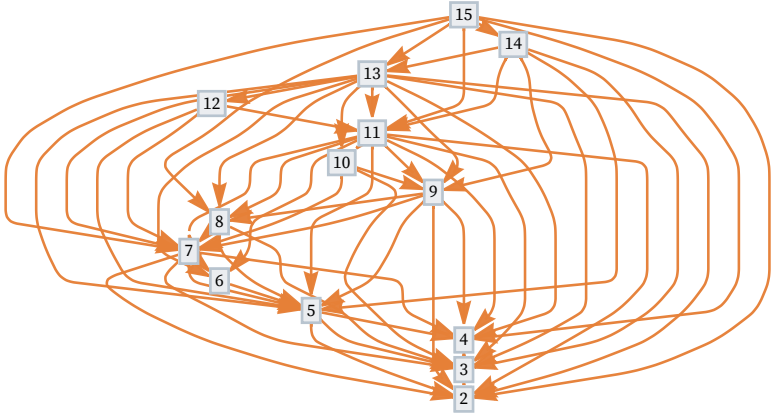

The transitive reduction of this is always trivial, however:

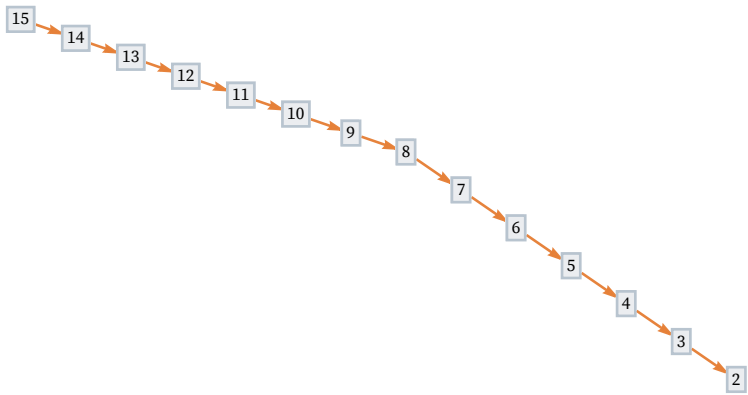

One general way to "probe" any function *f*[*n*] is to look at a multiway graph generated by the rule:

$n \longmapsto \{f[n], n + 1\}$

Here, for example, is the result for

$n \longmapsto \{\text{EulerPhi}[n], n + 1\}$

starting with $n = 1$:

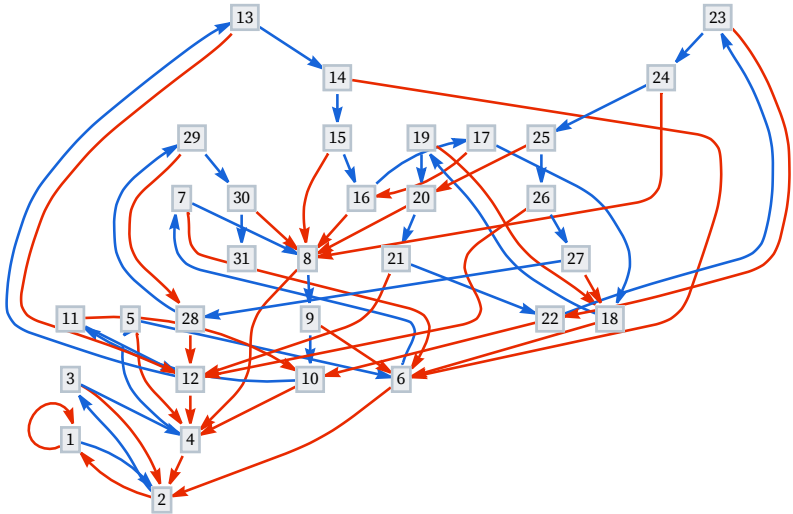

Once again, the transitive reduction is very simple:

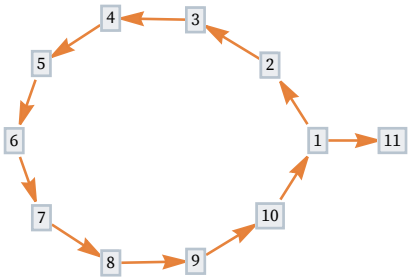

As another example, we can look at:

$n \longmapsto \{\text{PrimePi}[n], n + 1\}$

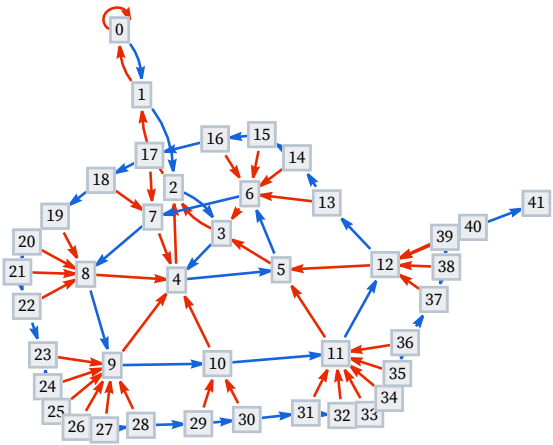

where each "efflorescence" corresponds to a prime gap:

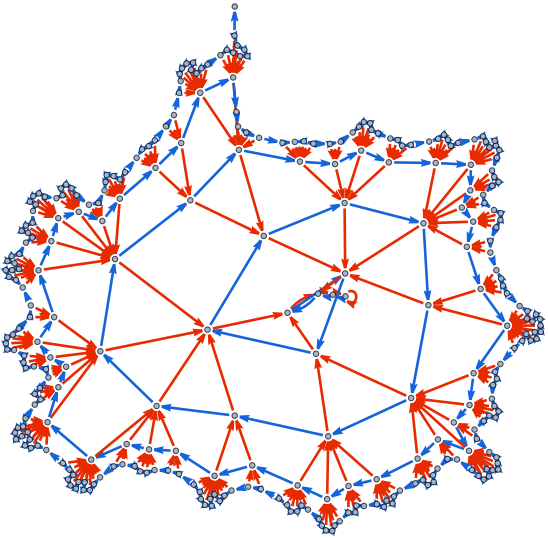

As a final example we can consider the digit-reversal function:

$n \longmapsto \{\text{IntegerReverse}[n, 2], n + 1\}$

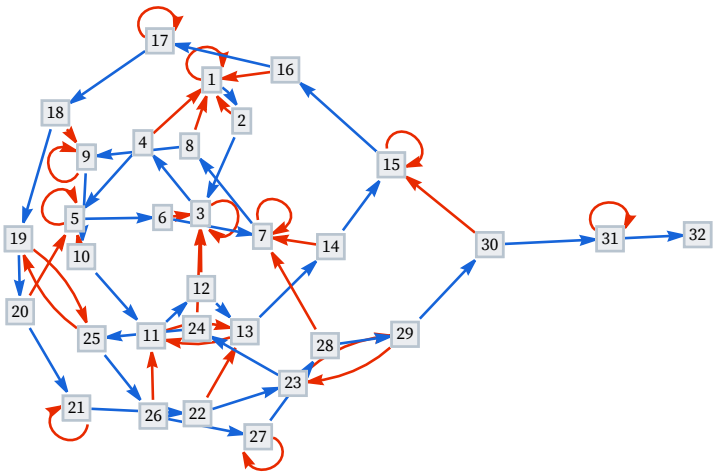

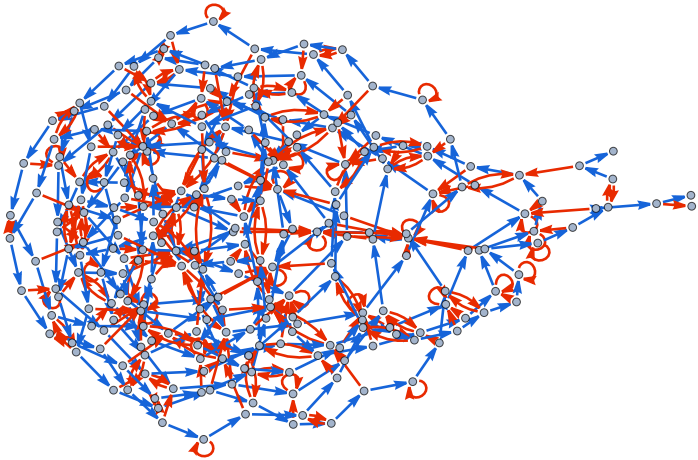

## Non-Integer Values

In almost everything we've discussed so far, we've been considering only integer values, both in our rules and our initial conditions. So what happens if we start a rule like

$n \longmapsto \{n + 1, n + 2\}$

with a non-integer value? Rather than taking a specific initial value, we can just use a symbolic value *x*—and it then turns out that the multiway graph is the same regardless of the value of *x*, integer or non-integer:

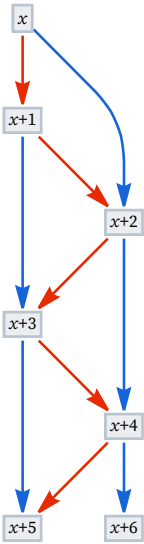

What if the rule contains non-integer values? In a case like

$$n \longmapsto \{n + a, n + b\}$$

the basic properties of addition ensure that the multiway graph will always have the same grid structure, regardless of *a*, *b* and the initial value *x*:

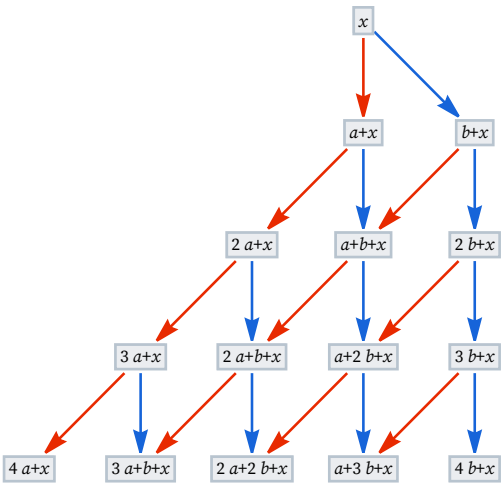

But in a case like

$$n \longmapsto \{n + a, b\,n\}$$

things are more complicated. For arbitrary symbolic $a$, $b$ and initial $x$, there are no relations that apply, and so the multiway graph is a pure tree:

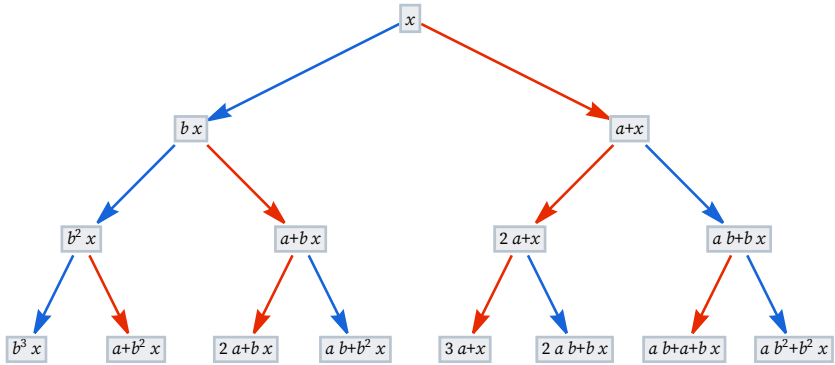

For a specific value of $b$, however, there are already relations, and a more complicated structure develops:

$n \longmapsto \{n + a, 2\,n\}$

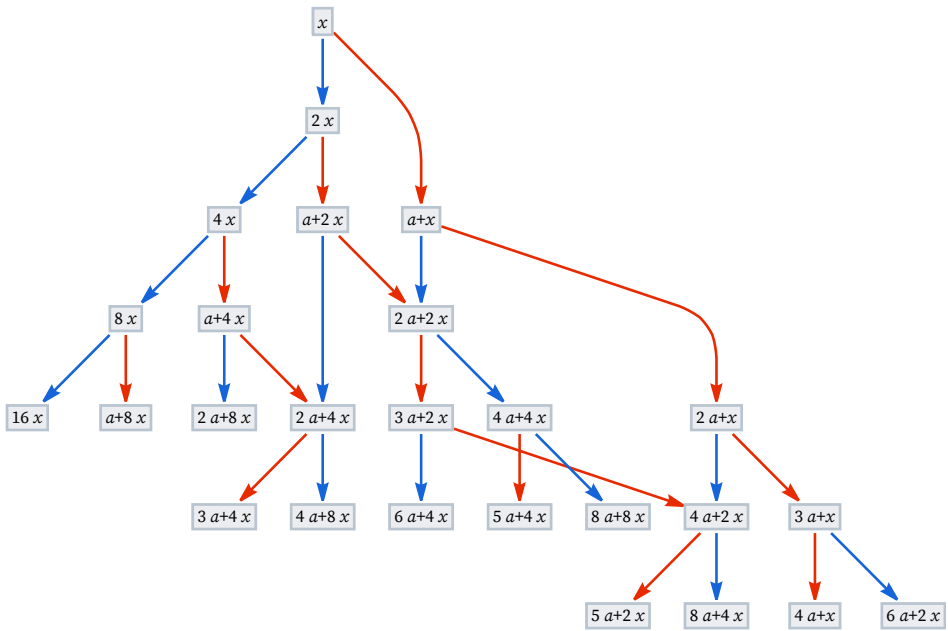

Continuing for more steps and removing loose ends we get

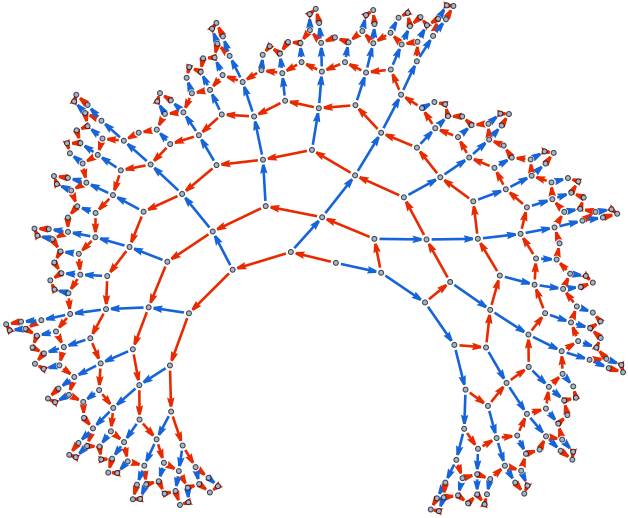

which is to be compared to the result from above for $a = 1$, $x = 1$:

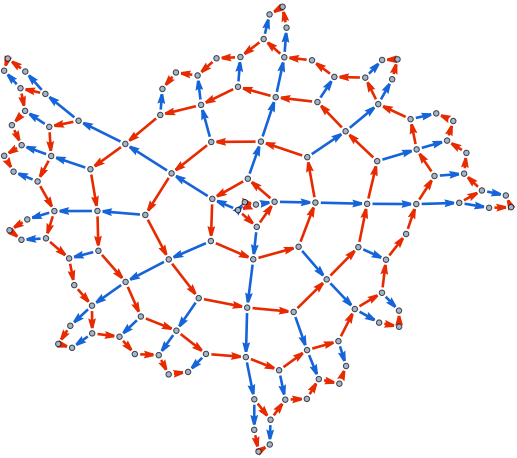

What happens if we choose a non-integer value of $b$, say:

$n \longmapsto \{n + 1, \ \sqrt{2} \ n\}$

We immediately see that there are "special relations" associated with $\sqrt{2}$ and its powers:

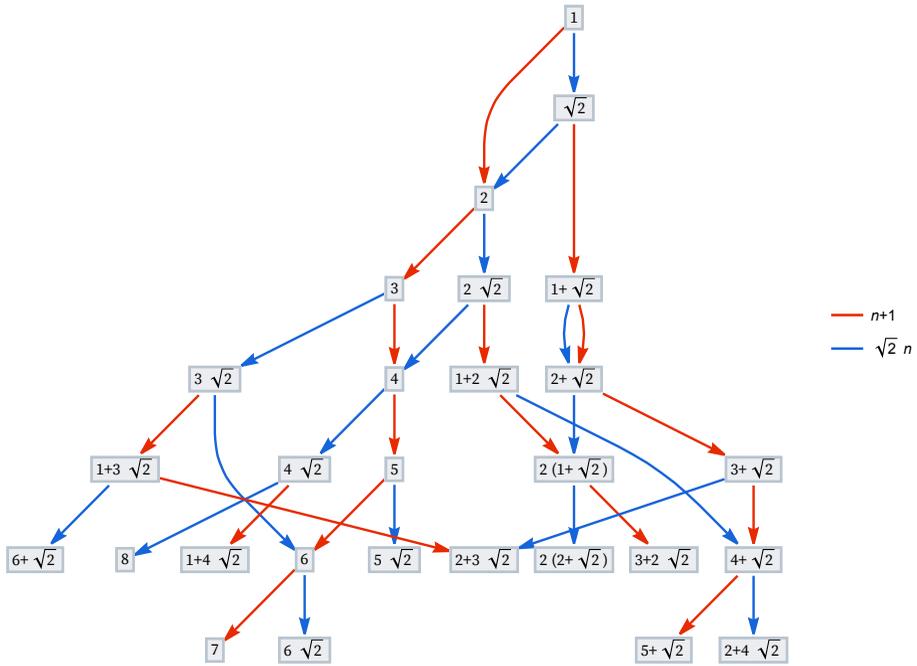

Continuing for longer we get the somewhat complex structure:

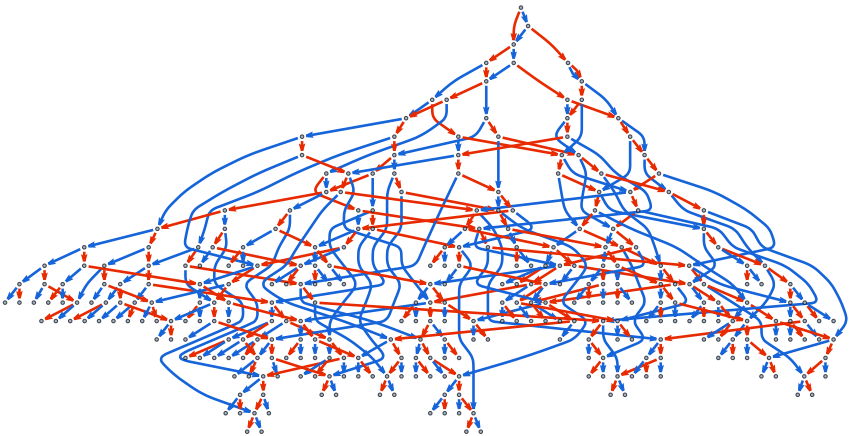

or in a different rendering with loose ends removed:

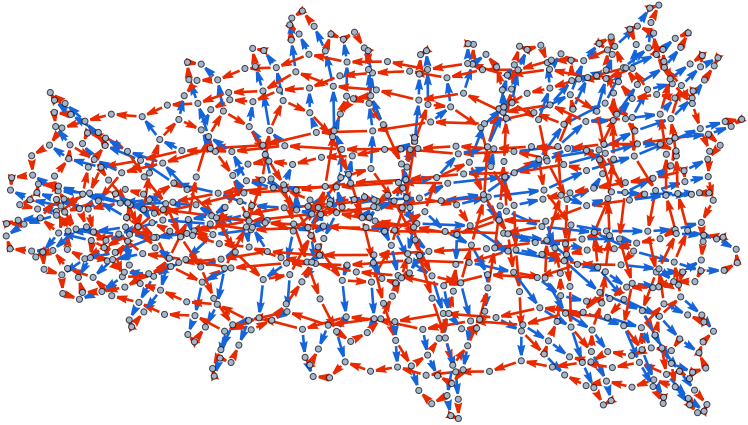

This structure is very dependent on the algebraic properties of $\sqrt{2}$. For a transcendental number like $\pi$ there are no "special relations", and the multiway graph will be a tree. For $\sqrt{3}$ we get

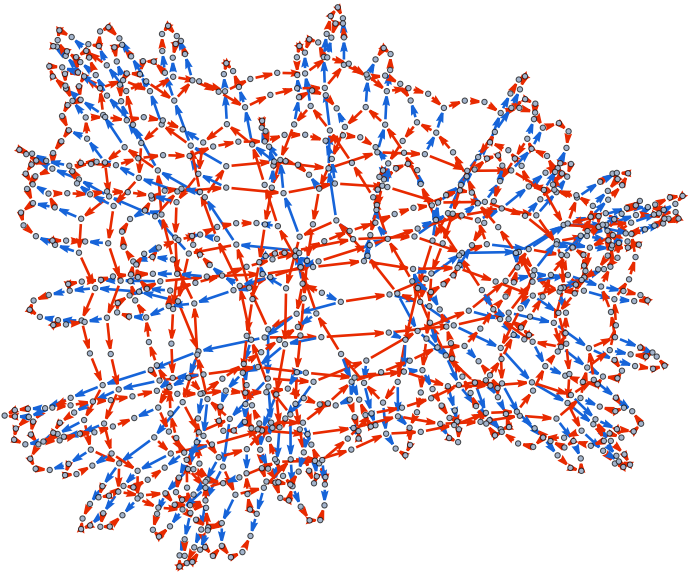

and for $1 + \sqrt{2}$ :

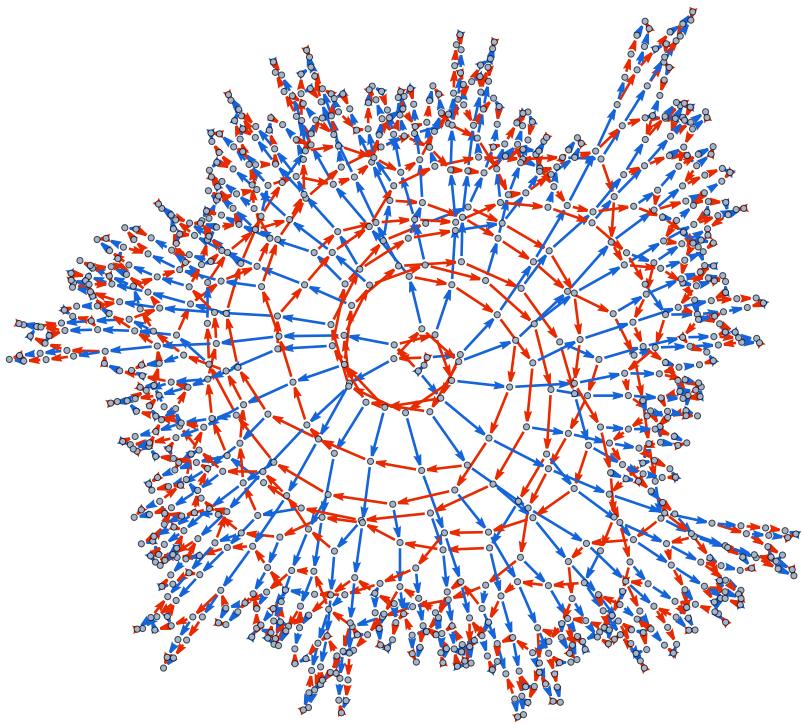

# Complex Numbers

There are many possible generalizations to consider. An immediate one is to complex integers.

For real numbers $n \mapsto \{a\,n, b\,n\}$ always generates a grid. But for example

$n \mapsto \{i\,n, 2\,n,\}$

instead generates

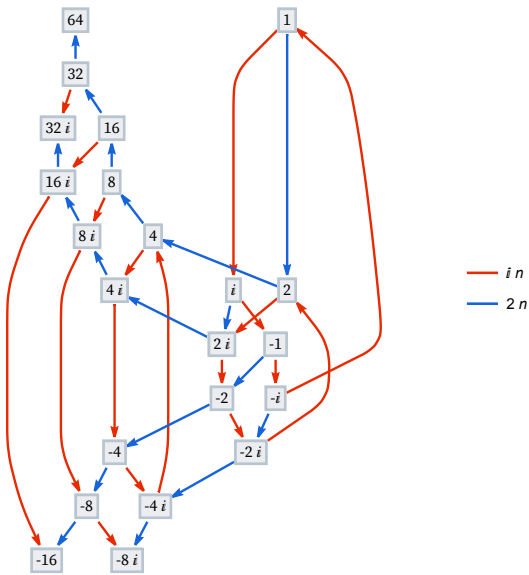

Continuing for longer, the graph becomes:

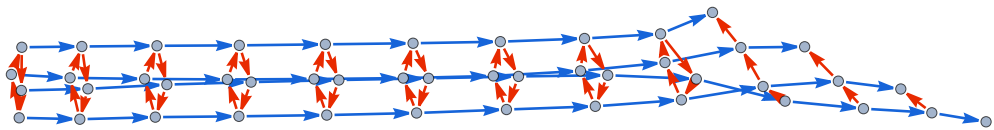

One feature of having values that are complex numbers is that these values themselves can be used to define coordinates to lay out the nodes of the multiway graph in the plane—giving in this case:

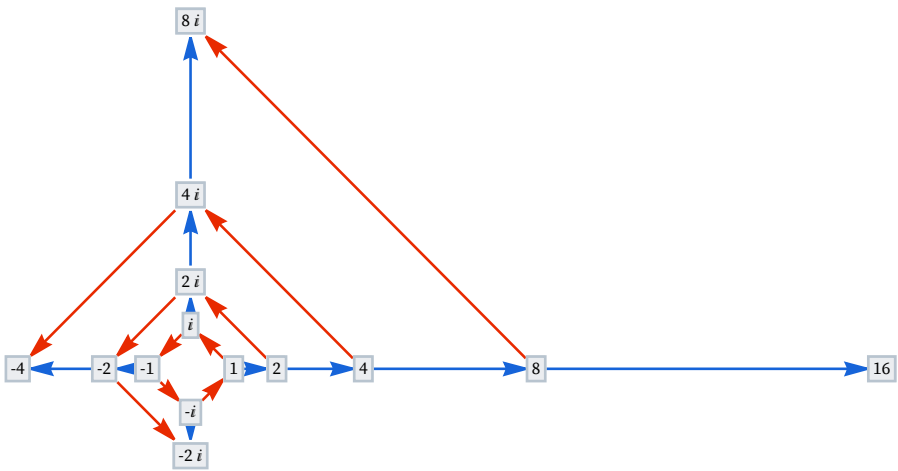

or after more steps:

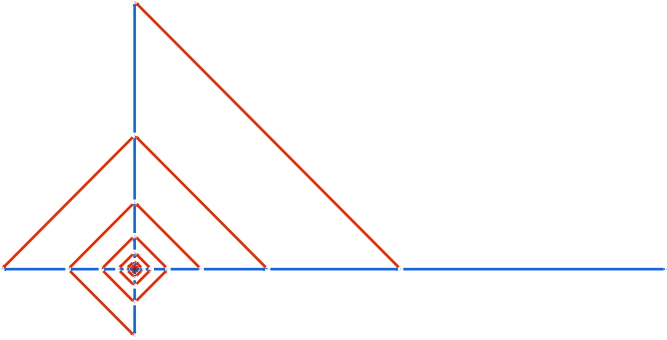

Similarly

$n \longmapsto \{(1 + i)\, n, \, (1 - i)\, n\}$

gives

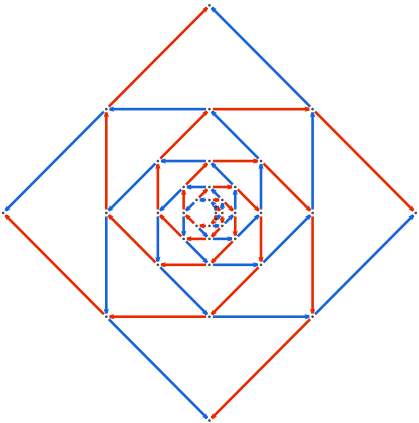

The non-branching rule

$n \longmapsto \{(1 + i)\, n\}$

yields

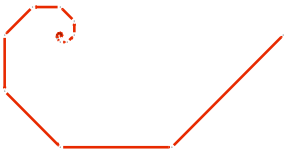

while

$n \longmapsto \{(3 + i)\, n, \, (2 - i)\, n\}$

gives

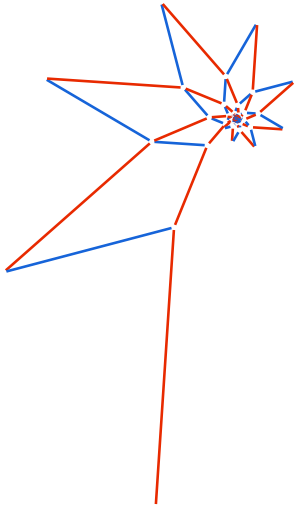

If we combine multiplication with addition, we get different forms—and we can make some interesting mathematical connections. Consider rules of the form

$$n \longmapsto \{1 + c\,n,\ 1 + \mathsf{Conjugate}[c]\,n\}$$

where $c$ is some complex number. I considered such rules in *A New Kind of Science* as a practical model of plant growth (though already then I recognized their connection to multiway systems). If we look at the case

$$c = \frac{1}{2} - \frac{i}{4}$$

the multiway graph is structurally just a tree:

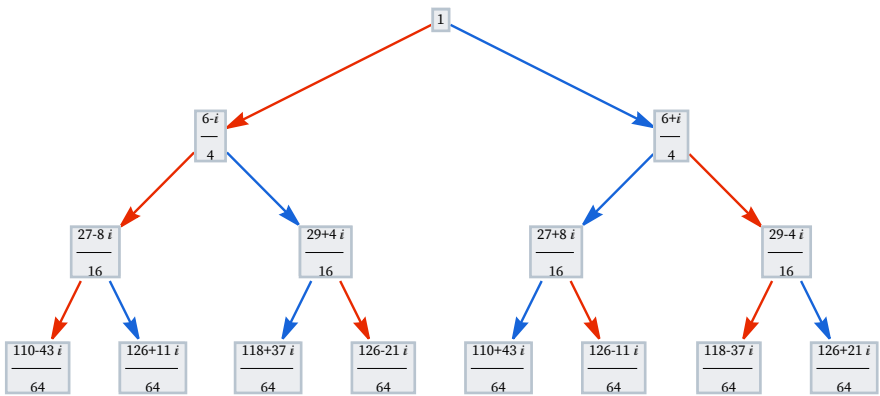

But if we plot nodes at the positions in the complex plane corresponding to their values we get:

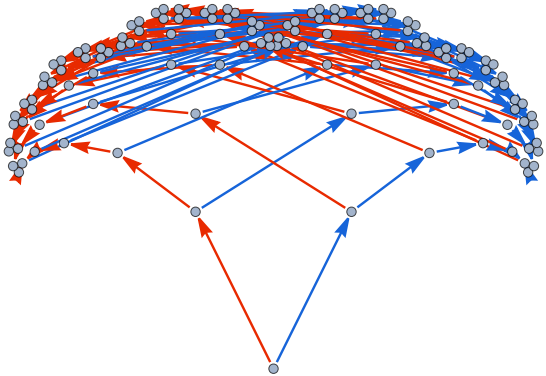

Continuing this, and deemphasizing the "multiway edges" we see a characteristic "fractal-like" pattern:

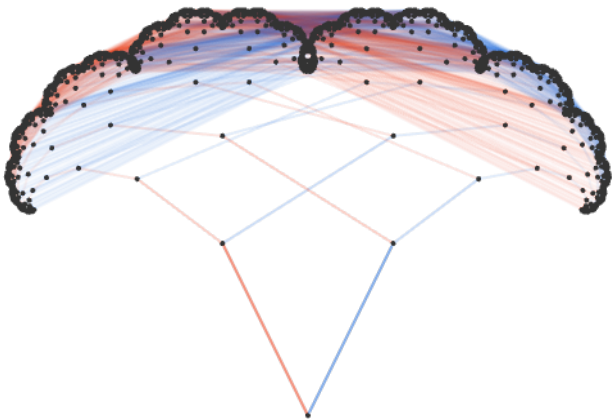

Note that this is in some sense dual to the typical "line segment iteration" nested construction:

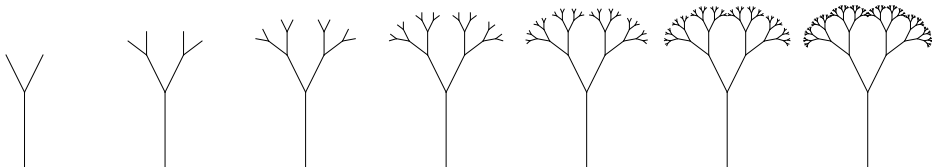

Adding a third "real" branch

$$n \mapsto \{1 - \frac{i\,n}{2}, 1 + \frac{i\,n}{2}, 1 + \frac{n}{2}\}$$

we get

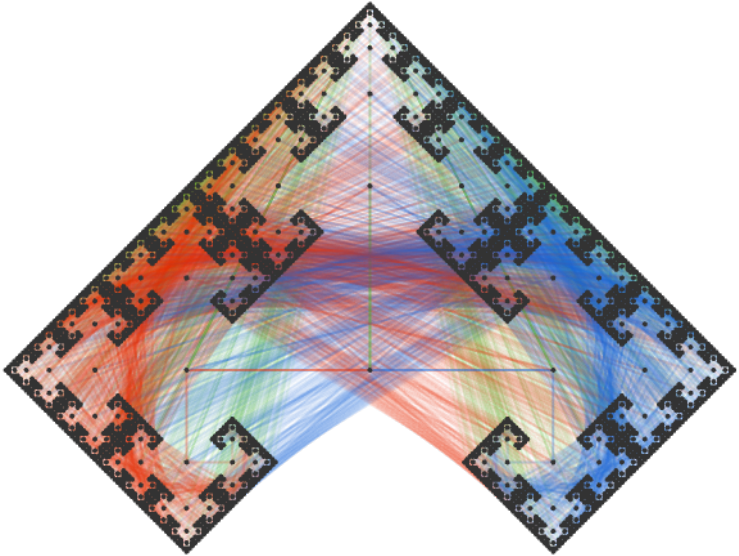

And with

$$n \mapsto \{1 + \frac{1}{4}(-1 - i\sqrt{3})\,n,\ 1 + \frac{1}{4}(-1 + i\sqrt{3})\,n,\ 1 + \frac{n}{2}\}$$

the result builds up to a typical Sierpinski pattern:

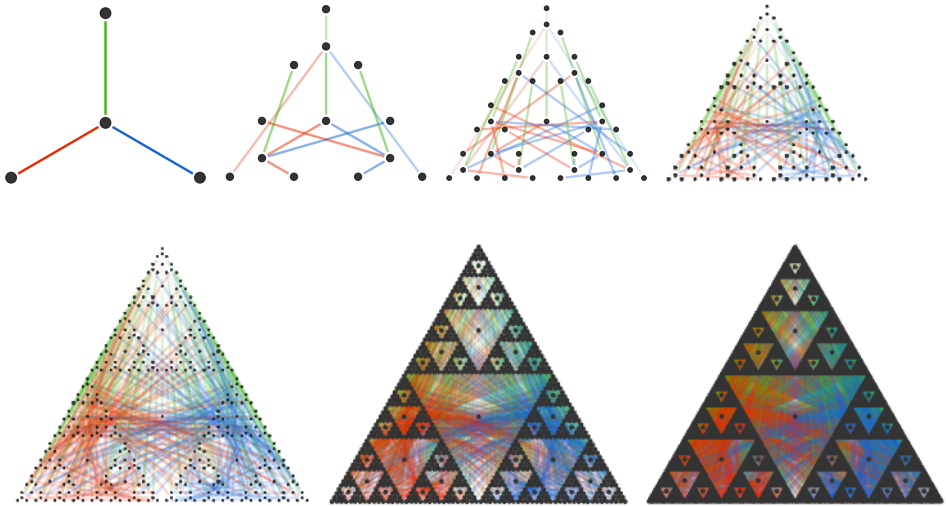

These pictures suggest that at least in the limit of an infinite number of steps there will be all sorts of merging between branches. And indeed it is fairly straightforward to prove this. But what about after, say, $t$ steps?

The result from each branch for the rule

$n \mapsto \{1 + (a + b\,\iota)\,n,\ 1 + (a - b\,\iota)\,n\}$

is a polynomial such as

$1 + (1 + (a - i\,b)\,(1 + a + i\,b))\,(a + i\,b)$

or

$1 + a + a^2 + a^3 + i\,b + i\,a^2\,b + b^2 + a\,b^2 + i\,b^3$

So now the question of merging becomes a question of finding solutions to equations which equate the polynomials associated with different possible branches. The simplest nontrivial case equates branch {1, 1} with branch {2, 2}, yielding the equation:

$1 + a + i\,b = 1 + 2\,a - 2\,i\,b$

with solution

$$a = -\frac{1}{2}$$

We can see this merging in action with the rule:

$n \mapsto \{1 + (-\frac{1}{2} + i)\,n,\ 1 + (-\frac{1}{2} - i)\,n\}$

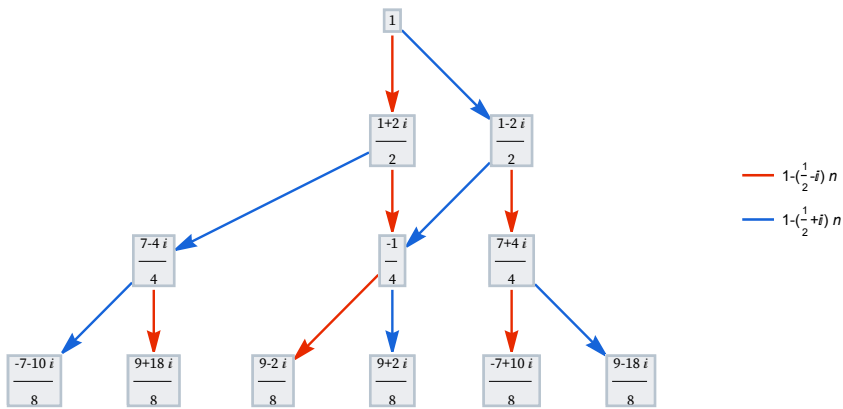

The core of what it generates is the repetitive structure:

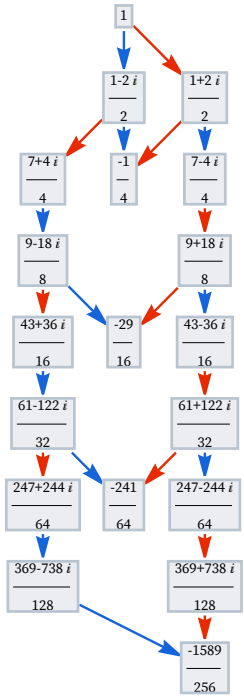

A few additional results are (where the decimals are algebraic numbers of degree 6, and *a* is a real number):

| | |
|---|---|
| {1, 22} | $-1 + i$ |
| {1, 112} | $-1 - i$ |
| {1, 222} | $-0.319448 - 1.29122\,i$ |
| {11, 112} | $-\dfrac{2}{3} + \dfrac{i\,\sqrt{2}}{3}$ |
| {11, 222} | $i\,\sqrt{2}$ |
| {12, 111} | $-0.726699 - 1.45975\,i$ |
| {111, 112} | $-\dfrac{1}{2} + \dfrac{i\,\sqrt{3}}{2}$ |
| {111, 222} | $a + i\,\sqrt{1 + 2\,a + 3\,a^2}$ |
| {112, 221} | $a + i\,\sqrt{1 - a^2}$   if $-1 < a < 1$ |

In a case like

$$n \mapsto \{1 + i\,\sqrt{2}\,n,\ 1 - i\,\sqrt{2}\,n\}$$

there is an "early merger"

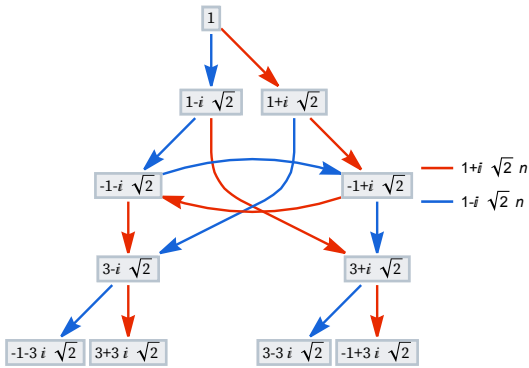

but then the system just generates a tree:

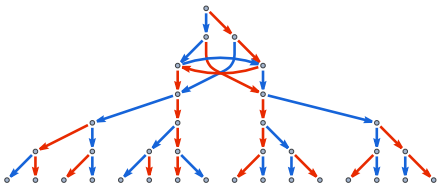

The family of rules of the form

$$n \mapsto \{1 + (a + i\sqrt{1-a^2})\, n,\ 1 + (a - i\sqrt{1-a^2})\, n\}$$

shows more elaborate behavior. For $a = \frac{1}{2}$ we get:

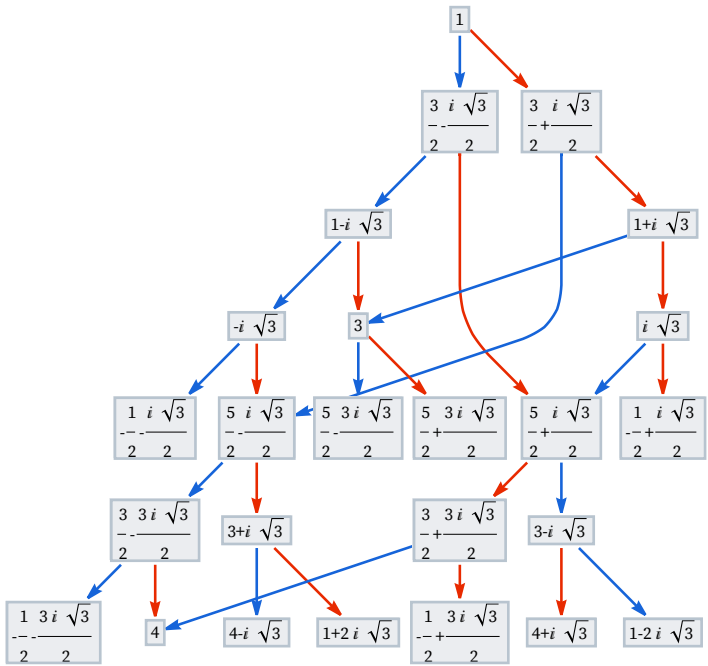

Continuing for more steps this becomes:

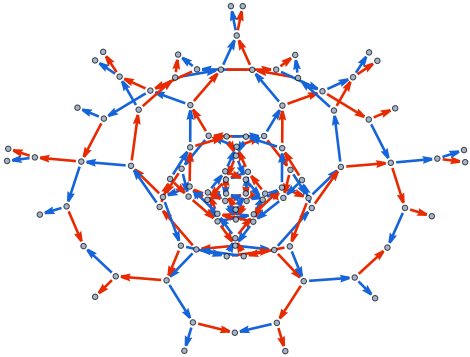

For $a = \frac{1}{4}$ we get instead:

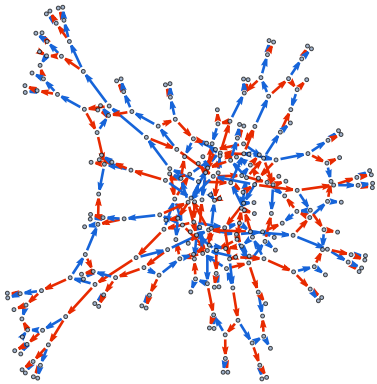

If we look at the actual distribution of values obtained by such rules we find for example:

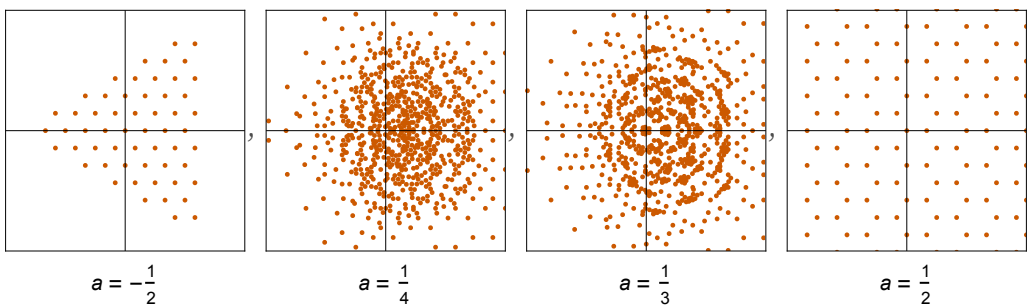

If we go beyond multiway systems with pure "$1 + c\,n$" rules we soon get results very similar to ones we've seen in previous sections. For example

$$n \longmapsto \{2\,|\,n, 1 + n\}$$

gives multiway graph (after removing loose ends)

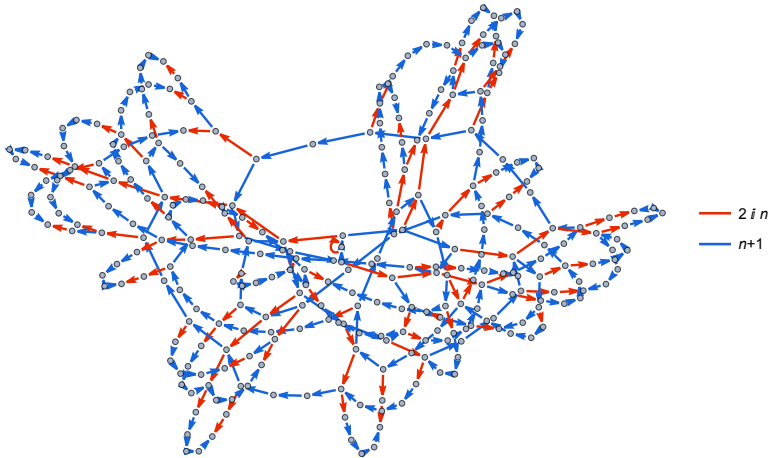

Placing nodes according to their numerical values this then has the form:

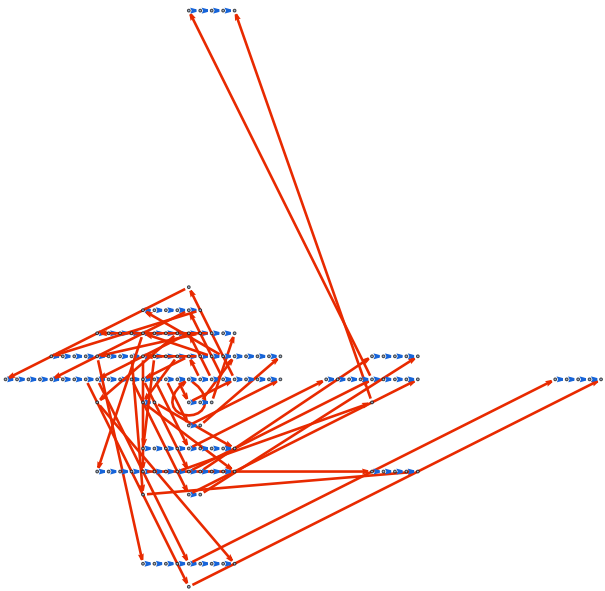

# Collections of Numbers, and Causal Graphs

In studying multiway systems based on complex numbers we're effectively considering a special case of multiway systems based on collections of numbers. If the complex-number rules are linear, then what we have are iterated affine maps—that form the basis for what I've called geometric substitution systems.

As a slightly more general case we can consider multiway systems in which we take pairs of numbers $v$ and apply the rule

$$v \mapsto \{a.v, b.v\}$$

where now $a$ and $b$ are 2×2 matrices. If both matrices are the form $\{\{r, -s\}, \{s, r\}\}$ then this is equivalent to the case of complex numbers. But we can also for example consider a rule like

$$v \mapsto \{\{\{1, 0\}, \{1, 1\}\}.v, \{\{1, 1\}, \{1, 1\}\}.v\}$$

which yields

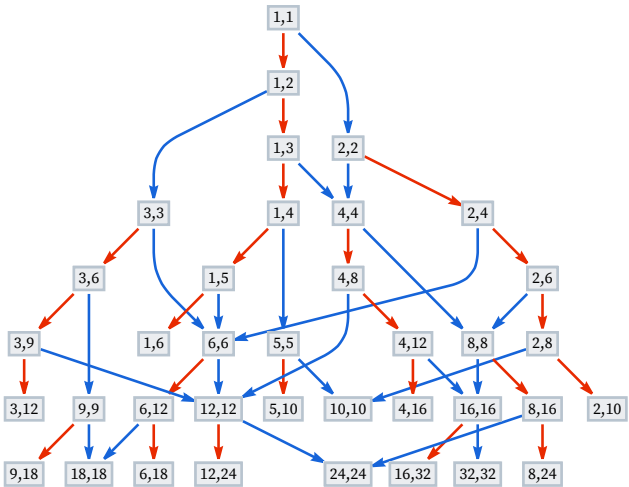

or after more steps and in a different rendering:

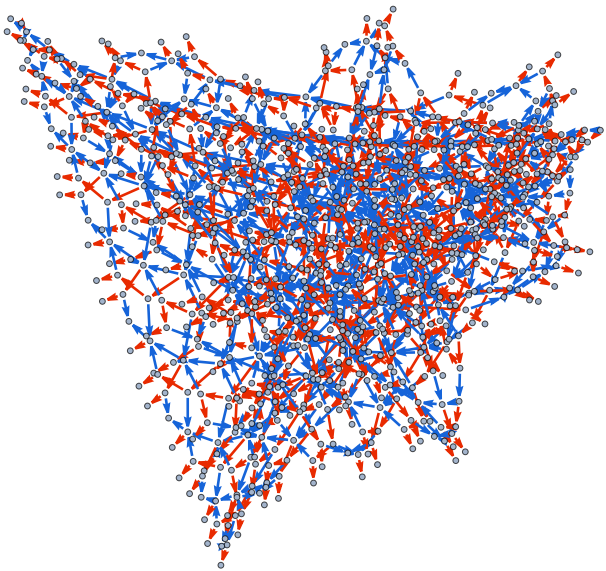

Laying this out in 2D using the actual pairs of numbers as coordinates, this becomes:

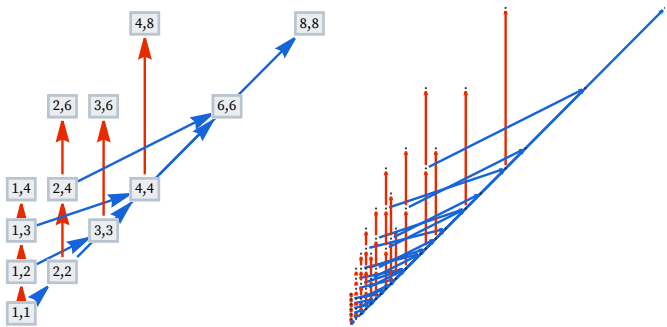

Here are samples of typical behavior with 2×2 0, 1 matrices:

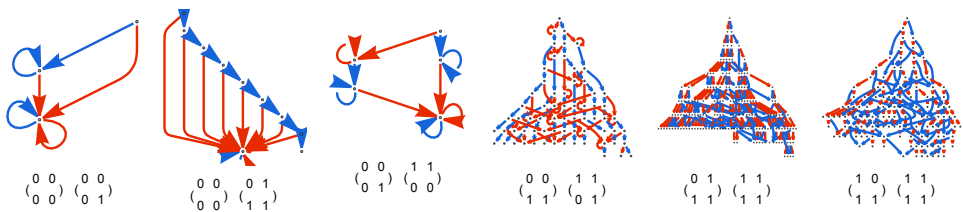

Beyond pure matrix multiplication, we can also consider a rule that adds constant vectors, as in:

$$v \mapsto \{a.v + c, b.v + d\}$$

We can also think in a more "elementwise" way, constructing for example simple rules such as

$$\{x, y\} \rightarrow \{\{x + 2, y + 1\}, \{y, x + 1\}\}$$

This generates the multiway graph:

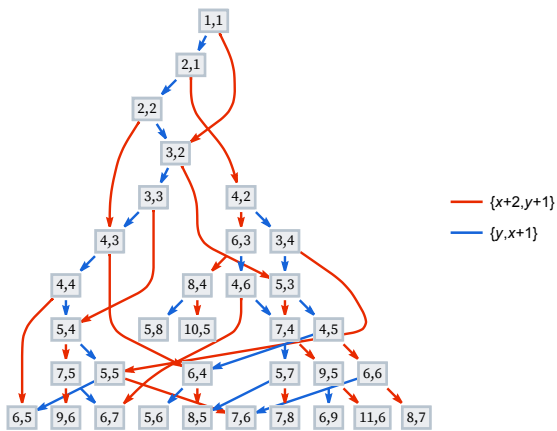

Continuing for longer and removing loose ends yields:

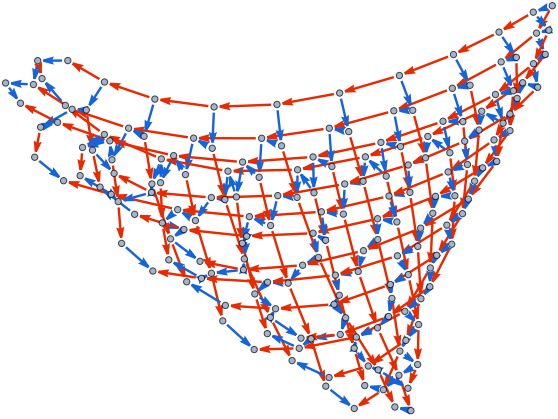

Using values as coordinates then gives:

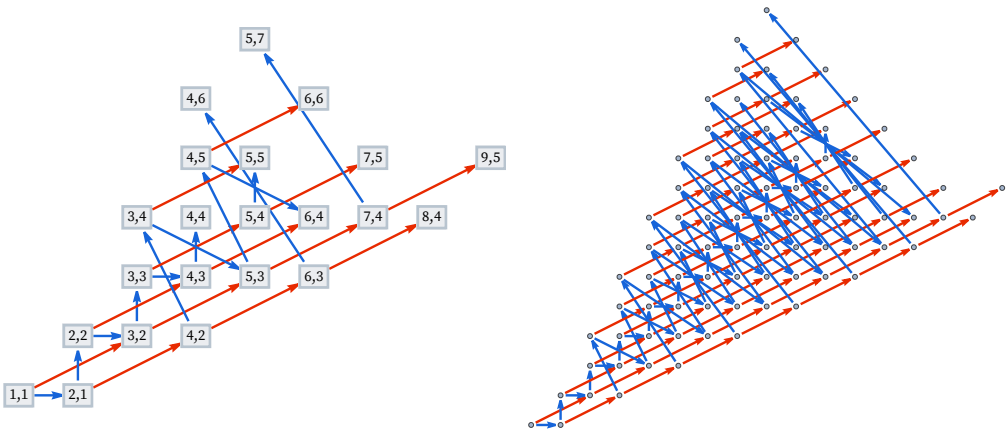

In our Physics Project and other applications of multicomputation, we often discuss causal graphs, that track the causal relationships between updating events. So why is it that these haven't come up in our discussion of multiway systems based on numbers? The basic reason is that when our states are individual numbers, there's no reason to separately track updating events and transformations of states because these are exactly the same—because every time a state (i.e. a number) is transformed the number as a whole is "consumed" and new numbers are produced. Or, in other words, the flow of "data" is the same as the flow of "causal information"—so that if we did record events, there'd just be one on each edge of the multiway graph.

But the story is different as soon as our states don't just contain individual "atomic" things, like single numbers. Because then an updating event can affect just part of a state—and asking what causal relationships there may be between events becomes something separate from asking about the transformation of whole states.

With a rule of the form, say,

$\{x, y\} \to \{f_1[x], \; f_2[x], \; g_1[y], \; g_2[x]\}$

things are still fairly trivial. Yes, there are separate "$x$" and "$y$" events. But they don't mix, so we'll just get two independent causal graphs. Things can be less trivial in a case like the one above, of the form:

$\{x, y\} \to \{f[x, y], \; g[x, y]\}$

But now there is a different problem. Let's say that the rule transforms $\{x, y\}$ to $\{y + 1, x + 1\}$. How should we decompose that into "elementary events"? We could say there's one event that swaps $x$ and $y$, and others that add 1. Or something different. It's hard to know.

So why haven't we encountered this kind of problem in other multicomputational systems, say in hypergraph rewriting systems or string substitution systems? The point is that in these systems the underlying elements always have a certain unique identity, which allows their "flow" to be traced. In our Physics Project, for example, each hypergraph updating event that occurs affects certain particular "atoms of space" (that we can think of as being labeled by unique identifiers)—and so we can readily trace how the effects of different events are related. Similarly, in a string substitution system, we can trace which characters at which positions in the string were affected by a given event, and we can then trace which new characters at which new positions these affect.

But in a system based on numbers this tracing of "unique elements" doesn't really apply. We might think of 3 as being 1+1+1. But there's nothing that uniquely tags these 1s, and allows us to trace how they affect 1s that might make up other numbers. In a sense, the whole point of numbers is to abstract away from the labeling of individual objects—and just ask the aggregate question of "how many" there are. So in effect the "packaging" of information into numbers can be thought of as "washing out" causal relationships.

When we give a rule based on numbers what it primarily does is to specify transformations for values. But it's perfectly possible to add an ancillary "causal rule", that, for example, can define which elements in an "input" list of numbers should be thought of as being "used as the inputs" to produce particular numbers in an output list of numbers.

There's another subtlety here, though. The point of a multiway graph is to represent all possible different histories for a system, corresponding to all possible sequences of transformations for states. A particular history corresponds to a particular path in the multiway graph. And if—as in a multiway system based on single numbers—each step in this path is associated with a single, specific event, then the causal graph associated with a particular history will always be trivial.

But in something like a hypergraph- or string-based system there's usually a nontrivial causal graph even for a single path of history. And the reason is that each transformation between states can involve multiple events—acting on different parts of the state—and there can be nontrivial causal relationships between these events "mediated" by shared elements in the state.

One can think of the resulting causal graph as representing causal relationships in "space-time". Successive events define the passage of time. And the layout of different elements in each state can be thought of as defining something like space. But in a multiway system based on single numbers, there isn't a natural notion of space associated with each state, because the states are just single numbers which "don't have enough structure" to correspond to something like space.

If we're dealing with collections of numbers, there's more possibility of "having something like space". But it's easiest to imagine this when one's dealing with very large collections of numbers, and when the "locations" of the numbers are more important than their values—at which point the fact that they're numbers (rather than, say, characters in a string) doesn't make much difference.

But in a multiway system one's dealing with multiple paths of history, not just one. And one can then start asking about causal relationships not just within a single path of history, but across different paths: a multiway causal graph. And that's the kind of causal graph we'll readily construct for a multiway system based on numbers. For a system based on strings or hypergraphs there's a certain wastefulness to starting with a standard multiway graph of transformations between states. Because if one looks at all possible states, there's typically a lot of repetition between the "context" of different updating events.

And so an alternative approach is to look just as the "tokens" that are involved in each event: hyperedges in a hypergraph, or runs of characters in a string. So how does it work for a multiway system based on numbers? For this we have to again think about how our states are decomposed for purposes of events, or, in other words, what the "tokens" in them are. And for multiway systems based on single numbers, the natural thing is just to consider each number as a token.

For collections of numbers, it's less obvious how things should work. And one possibility is to treat each number in the collection as a separate token, and perhaps to ignore any ordering or placement in the collection. We could then end up with a "multi-token" rule like

$\{n\_, m\_\} \mapsto \{n + m, n - m\}$

whose behavior we can represent with a token-event graph:

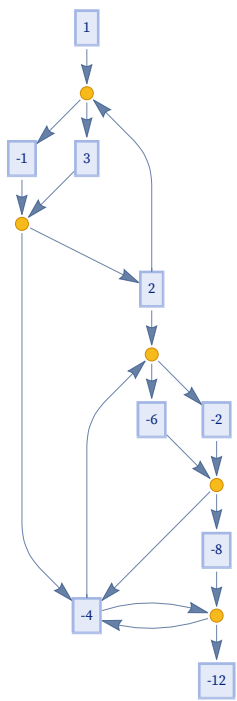

But given this, there is then the issue of deciding how collections of tokens should be thought of as aggregated into states. And in general multi-token numerical multiway systems represent a whole separate domain of exploration from what we have considered here.

A basic point, however, is that while our investigations of things like hypergraph and string systems have usually had a substantial "spatial component", our investigation of multiway systems based on numbers tends to be "more branchial", and very much centered around the relationships between different branches of history. This does not mean that there is nothing "geometrical" about what is going on. And in fact we fully expect that in an appropriate limit branchial space will indeed have a geometrical structure—and we have even seen examples of this here. It is just that that geometrical structure is—in the language of physics—about the space of quantum states, not about physical space. So this means that our intuition about ordinary physical space won't necessarily apply. But the important point is that by studying multiway systems based on numbers we can now hope to sharpen our understanding and intuition about things like quantum mechanics.

## Much More to Explore…

The basic setup for multiway systems based on numbers is very simple. But what we've seen here is that—just like for so many other kinds of systems in the computational universe—the behavior of multiway systems based on numbers can be far from simple.

In many ways, what's here just scratches the surface of multiway systems based on numbers. There is much more to explore, in many different directions. There are many additional connections to traditional mathematics (and notably number theory) to be made. There are also questions about the geometrical structures that can be generated, and their mathematical characterization.

In the general study of multicomputational systems, branchial—and causal—graphs are important. But here we have barely begun to consider them. A particularly important issue that we haven't addressed at all is that of alternative possible foliations. In general it has been difficult to characterize these. But it seems possible that in multiway systems based on numbers these may be amenable to investigation with some kind of mathematical techniques. In addition, for things like our Physics Project questions about the coordinatization of branchial space are of great significance—and the "natural coordinatizability" of numbers makes multiway systems based on numbers potentially an attractive place to study these kinds of questions.

Here we've considered only ordinary multiway systems, in which the rules always transform one object into several. It's also perfectly possible to study more general multicomputational systems in which the rules can "consume" multiple objects—and this is particularly straightforward to set up in the case of numbers.

Here we've mostly looked at multiway systems whose states are individual integers. But we can consider other kinds of numbers and collections of numbers. We can also imagine generalizing to other kinds of mathematical objects. These could be algebraic constructs (such a polynomials) based on ordinary real or complex numbers. But they could also, for example, be objects from universal algebra. The basic setup for multiway systems—involving repeatedly applying functions—can be thought of as equivalent to repeatedly multiplying by elements (say, generators) of a semigroup. Without any relations between these elements, the multiway graphs we'll get will always be trees. But if we add relations things can be more complicated.

Multiway systems based on semigroups are in a sense "lower level" than ones based on numbers. In something like arithmetic, one already has immediate knowledge of operations and equivalences between objects. But in a semigroup, these all have to be built up. Of course, if one goes beyond integers, equivalences can be difficult to determine even between numbers (say different representations of radicals or, worse, transcendental numbers).

In their basic construction, multiway systems are fundamentally discrete—involving as they do discrete states, discrete branches, and discrete notions like merging. But in our Physics Project and other applications of the multicomputational paradigm it's often of interest to think about "continuum limits" of multiway systems. And given that real numbers provide the quintessential example of a continuum one might suppose that by somehow looking at multiway systems based on real numbers one could understand their continuum limit.

But it's not so simple. Yes, one can imagine allowing a whole "real parameter's worth" of outputs from the multiway rule. But the issue is how to "knit these together" from one step to the next. The situation is somewhat similar to what happens when one looks at ensembles of random walks, or stochastic partial differential equations. But with multiway systems things are both cleaner and more general. The closest analogy is probably to path integrals of the kind considered in quantum mechanics. And in a sense this is not surprising, because it is precisely the appearance of multiway systems in our Physics Project that seems to lead to quantum mechanics—and in a "continuum limit" to the path integral there.

It's not clear just how multiway systems are best generalized to the continuum case. But multiway systems based on numbers seem to provide a potentially promising bridge to existing mathematical investigations of the continuum—and I think have a good chance of revealing some elegant and powerful mathematics.

I first looked at multiway systems based on numbers back in the early 1990s, and I always meant to come back and look at them further. But what we've found here is that they're richer and more interesting than I ever imagined. And particularly from what we've now seen I expect them to have a very bright future, and for all sorts of important science and mathematics to connect to them, and flow from them.

## Thanks

I worked on what's described here during two distinct periods: May 2020 and September 2021. I thank for help of various kinds Tali Beynon, José Manuel Rodríguez Caballero, Bernat Espigule-Pons, Jonathan Gorard, Eliza Morton, Nik Murzin, Ed Pegg and Joseph Stocke—as well as my weekly virtual high-school "Computational Adventures" group.

---

*Cite as:* S. Wolfram (2021), "Multicomputation with Numbers: The Case of Simple Multiway Systems". https://www.wolframphysics.org/bulletins/2021/10/multicomputation-with-numbers-the-case-of-simple-multiway-systems.